\documentclass[a4paper,10pt]{article}

\usepackage[english]{babel}
\usepackage[T1]{fontenc}
\usepackage[dvips]{graphicx}
\usepackage{amsmath, amsfonts,amsthm, mathrsfs, amssymb}
\usepackage{verbatim}

\usepackage{enumerate}
\usepackage[latin1]{inputenc}
\usepackage{geometry}

\def\Im{{\rm Im}}
\def\Re{{\rm Re}}
\def\vars#1{\varsigma^{(#1)}}
\def\idjd{\left|i^d-j^d\right|}
\def\sidjd{ \inter{i^d-j^d}}
\def\ak{1+|k|^\tau}
\def\vsi{\varsigma}

\def\ijkno{\left|i-j\right|+|k|\not=0}
\def\inte#1{\lfloor#1\rfloor}
\def\inter#1{\langle #1\rangle}

\def\vf{\varphi}

\def\eX{e^{-\im\epsilon X}}
\def\meX{e^{\im\epsilon X}}

\def\bee{{\bf e}}
\def\be#1{\bee_{#1}}

\def\ir{{\rm i}}

\def\Lip{\cL ip}
\def\diag{{\rm diag}}
\def\normaL#1{\norma{#1}^{\cL}}

\def\norma#1{\left\|#1\right\|}
\def\sleq{\preceq}
\def\sgeq{\succeq}

\def\hz{h_0}
\def\widebar#1{\overline{#1}}

\newcommand{\C}{{\mathbb C}}

\newcommand{\N}{{\mathbb N}}

\newcommand{\R}{{\mathbb R}}

\newcommand{\T}{{\mathbb T}}
\newcommand{\Z}{{\mathbb Z}}

\newcommand{\cB}{{\mathcal B}}

\newcommand{\cH}{{\mathcal H}}

\newcommand{\cK}{{\mathcal K}}
\newcommand{\cL}{{\mathcal L}}

\newcommand{\cO}{{\mathcal O}}

\newcommand{\cQ}{{\mathcal Q}}
\newcommand{\cR}{{\mathcal R}}
\newcommand{\cS}{{\mathcal S}}

\newcommand{\norm}[1]{\| #1 \|}

\newcommand{\im}{{\rm i}}

\def\uno{{\bf 1}}

\numberwithin{equation}{section}
\newtheorem{theorem}{Theorem}[section]
\newtheorem{example}[theorem]{Example}
\newtheorem{lemma}[theorem]{Lemma}
\newtheorem{corollary}[theorem]{Corollary}

\newtheorem{definition}[theorem]{Definition}
\newtheorem{remark}[theorem]{Remark}

\title{Reducibility of 1-d Schr\"odinger equation with time quasiperiodic
  unbounded perturbations, I}

\author{D. Bambusi\footnote{Dipartimento di Matematica, Universit\`a degli Studi di Milano, Via Saldini 50, I-20133
Milano. \newline
 \textit{Email: } \texttt{dario.bambusi@unimi.it}},
}

\begin{document}

\maketitle

\begin{abstract}
We study the Schr\"odinger equation on $\R$ with a polynomial
potential behaving as $x^{2l}$ at infinity, $1\leq l\in\N$ and with a
small time quasiperiodic perturbation. We prove that if the symbol of
the perturbation grows at most like $(\xi^2+x^{2l})^{\beta/(2l)}$,
with $\beta<l+1$, then the system is reducible. Some extensions including
cases with $\beta=2l$ are also proved. The result implies boundedness
of Sobolev norms. The proof is based on pseudodifferential calculus
and KAM theory.
\end{abstract}

\section{Introduction}

In this paper we study the problem of reducibility of the time dependent
Schr\"odinger equation
\begin{align}
\label{schro}
\ir \dot\psi=H_\epsilon(\omega t)\psi\ , \ x\in\R
\\
\label{H}
H_\epsilon(\omega t):=-\partial_{xx}+V(x)+\epsilon W(x,-\ir \partial_x,\omega t)
\end{align}
where $V$ is a polynomial potential of degree $2l$, with $l\geq 1$,
and $\T^n\ni\phi\mapsto W(x,\xi,\phi)$ is a $C^\infty$ map from $\T^n$
to a space of symbols growing at infinity at most like
$(\xi^2+x^{2l})^{\beta/2l}$. We emphasize that the harmonic potential
$l=1$ is included. 

We will prove that, if $\beta<l+1$, then, for sufficiently small
$\epsilon$, and for $\omega$ belonging to a set of large measure,
there exists a unitary transformation which conjugates
Eq. \eqref{schro} to a time independent equation; the transformation
depends on time in a smooth quasiperiodic way. We also deduce
boundedness of the Sobolev norms and pure point spectrum of the
Floquet operator. In the case where the average of the symbol $W$ with
respect to the flow of the classical Hamiltonian system $\xi^2+V(x)$
vanishes, the result holds also for $\beta< (3l+1)/2$. Finally we
prove reducibility also in some cases with $\beta=2l$.

The main limitation of the paper is that the allowed perturbations are
of a quite particular type (it is the same as in \cite{HR82, HR82D}),
as an example, in the case
\begin{equation}
\label{nonp}
W=-\ir a_1(x,\omega t)\partial_x+a_0(x,\omega t)
\end{equation}
the functions $a_0$ and $a_1$ must be polynomials in $x$. On the
contrary the perturbation is allowed to grow at infinity (both in $x$
and in the Fourier variable $\xi$) much faster then in all the
preceeding papers.

There is quite an extensive literature on the problem of reducibility
of time dependent Schr\"odinger equation and the related problems of
growth of Sobolev norm and nature of the spectrum of the Floquet
operator. We recall first the works \cite{DS96,DSV02}, in which pure
point nature of the Floquet spectrum is obtained in the case in which
the growth of $V$ is superquadratic (and therefore the spectrum has
increasing gaps) and the perturbations is bounded and time
periodic. The first paper dealing with an unbounded time quasiperiodic
perturbation is \cite{BG01}. In \cite{BG01} we assumed that the
potential (not necessarily a polynomial) grows at infinity like
$x^{2l}$, with a real $l>1$ and the perturbation is bounded by
$1+|x|^\beta$ with $\beta<l-1$; reducibility in the limiting case
$\beta=l-1$ was obtained in \cite{LY10}. Concerning the case of
Harmonic potential we recall the pioneering work \cite{C87} in which
reducibility is obtained in case of a perturbation which is smoothing
and the works \cite{W08} and \cite{GT11} dealing with the case of a
bounded perturbation. The present paper is the first one in which
reducibility for an unbounded perturbation of the Harmonic oscillator
is obtained. We remark that the present result does not cover the
results of \cite{W08,GT11} since their perturbations are not in the
class of symbols we use here. The technique of the present paper can
be used also to obtain and improve \cite{W08,GT11}, but this requires
a quite heavy work and produces a bigger limitation for the allowed
range of $\beta$. For this reason it will be developped in a future
paper (paper II). 

We also recall the interesting counterexamples in \cite{GY00} and
\cite{D14}. In particular we remark that the class of perturbation
constructed in these papers is covered by the result of the present
paper\footnote{Actually in \cite{GY00} it is shown that the
  instability exhibited in their counterexample is stable under a
  class of further perturbations not covered by the present
  paper. This class will be covered in a paper II.}. The main
point is that in our case the frequencies fulfill a non-resonance
relation which is violated in \cite{GY00,D14}.

We recall that all the papers quoted above deal only with the one
dimensional case. The case of higher dimension is dealt with only in
the papers \cite{EK09} for the Schr\"odinger equation on $\T^d$ and in
\cite{GP16} for the case of the Harmonic oscillator. 

We remark that the problem of reducibility of linear equations is
considered to be the main step for the proof of KAM type results in
nonlinear PDEs, thus we
think that the result of the present paper could be useful in this
direction and in particular in order to construct quasiperiodic
motions of a soliton in external potentials (in the spirit of
\cite{FGJS04,BM16}).

The proof of the result of the present paper is based on a
generalization of the ideas developed by Baldi, Berti, Montalto
\cite{BBM14} (see also \cite{Mtesi,FP15,BM16a}) in
order to extend KAM theory to fully nonlinear equations, ideas which
in turn are a development of those introduced by Plotnikov and Toland
in \cite{PT01} in order to study the water wave problem (see also
\cite{IPT05}).  We recall that the idea is to proceed in two steps:
first one uses pseudodifferential calculus in order to regularize the
perturbation and then applies more or less standard KAM theory in
order to conclude the proof. Actually an intermediate step is also
required. This is due to the fact that, after the smoothing theorem
the system is reduced to a smoothing perturbation of a time independent
system, but the time independent system is not diagonal. So before
developing KAM theory one has to diagonalize such a time independent
system and to study its eigenvalues. 

The main novelty of the present paper is that we deal here with the
case of an equation on an unbounded domain, namely $\R$ so that a
second source of unboundedness is the growth at infinity of
symbols. In order to deal with the present case one has to develop in
a quite careful way the regularization procedure, which is based on
the strong connection existing between classical and quantum
perturbation theory \cite{GP87, BGP99}. The point is that, if one
considers the classical Hamiltonian of the system and tries to
eliminate order by order (in $\epsilon$) the time dependence through
classical normal form theory, then the quantization of the normalizing
transformation conjugates the quantum system to a time independent
system, {\it up to the quantum corrections}. But the quantum
corrections are usually smoother then the original operators, and
therefore one can expect the transformed quantum system to be a
smoother perturbation of a time independent system. It turns out that
this is the case. The framework (and the results) that
we use here is the one developed by Hellfert and Robert in
\cite{HR82}.

After the regularization step one can use more or less standard KAM
theory in order to reduce the regularized system to constant
coefficients.  However there is an additional difficulty, namely that
pseudodifferential calculus works well in the class of $C^\infty$
functions, while the simplest formulation of KAM theory is that
dealing with analytic functions. So one has to develop KAM theory in a
$C^\infty$ context. This is quite standard and indeed KAM theory is
developed in $C^\infty$ context e.g. in the paper \cite{BBM14},
however we are here in a slightely different situation, thus we
decided to insert in the paper also a proof of a KAM theorem
with finite smoothness developed following the presentation of
\cite{S04}. We point out that the method of \cite{S04} has already
been applied to the problem of reducibility, in a slightly different
context in \cite{YZ13}.

As anticipated above the main limitation of the present paper is that
the symbols we consider here are of a quite particular type. The
extension to more general symbols only fulfilling growth properties
will be the goal of the paper II. The main point in order to get
the extension is to introduce a different class of symbols; however,
on the one hand a quite hard technical work is needed in order to deal
with such a class, and on the other one we only get the result under the
strongest assumption $\beta<l$, which in particular rules out the case
$\beta=2l-1$ which is very interesting in order to deal with the case
of a soliton moving in an external potential.

\vskip 10 pt

The paper is organized as follows: In sect. \ref{main} we state the
results of the paper and give some examples and comments. The
subsequent sections contain the corresponding proofs. Precisely, in
Sect. \ref{algo} we introduce and give the main properties of the
unitary transformations generated by time dependent selfadjoint
operators. Such transformations will be used in the rest of the paper
first at level of symbols and subsequently directly at the level of
operators. In sect. \ref{egor} we prove the smoothing theorem. The
section is split into a few subsections. In particular, in Subsection
\ref{smoa} there is quite detailed description of the strategy used in
order to prove the smoothing theorem. In Sect. \ref{prepa} we
diagonalize the time independet part of the regularized system and
study its eigenvalues. In Sect. \ref{lemmaKAM} we prove the analytic
KAM theorem that constitutes the main step for the proof of finite
smoothness KAM theorem proved in sect. \ref{finite}. Finally, the
appendix contains some technical Lemmas. They are grouped in some
different sections according to the role they have in the main part of
the text.

\noindent{\it Acknowledgements}. This paper originated from a series
of discussions with quite a lot of people on the methods of
\cite{BBM14,Mtesi,FP15,BM16a} and on the possibility of extending them
to the case of the Schr\"odinger equation. In particular I warmly
thank P. Baldi, R. Montalto and M. Procesi who explained to me in a
quite detailed way their works. During the preparation of the present
work I benefit of many suggestions and discussions with A. Maspero and
D. Robert. In particular D. Robert pointed to my attention (and often
explained me) his papers in which the class of symbols that are used
here are extensively studied. I also thank B. Gr\'ebert for some
discussions on the Harmonic case that allowed me to fix some points of
the proof.

\section{Statement of the Main Result}\label{main}

Fix a positive integer $l\geq 1$ and define the weight
\begin{equation}
\label{lam}
\lambda(x,\xi):=\left(1+\xi^2+x^{2l}\right)^{1/2l}\ ,
\end{equation}
\begin{definition}
\label{Sm}
The space $S^{m}$ is the space of the symbols $g\in C^\infty(\R)$ such
that $\forall k_1,k_2\geq 0$ there exists $C_{k_1,k_2}$ with the
property that\footnote{This class of symbols coincides with the class
  introduced in \cite{HR82} and denoted by $S^m_{l,1}$. }
\begin{equation}
\label{sm}
\left|\partial^{k_1}_\xi\partial^{k_2}_{x}g(x,\xi)\right|\leq
C_{k_1,k_2} \left[\lambda(x,\xi)\right]^{m-k_1l-k_2}\ .
\end{equation}
The best constants $C_{k_1,k_2}$ such that \eqref{sm} hold form a
family of seminorms for that space $S^{m}$. 
\end{definition}

\begin{remark}
\label{smo}
All what we will do can be developed also for symbols with a finite,
but large, differentiability.
\end{remark}

In the following we will denote by
$\cS^{m}:=C^{\infty}(\T^n,S^{m})$ the space of
$C^{\infty}$ functions on $\T^n$ with values in $S^{m}$.

The frequencies $\omega$ will be assumed to vary in the set 
$$
\Omega:=[1,2]^{n}\ ,
$$
or in suitable closed subsets
$\widetilde\Omega$. 

To a symbol $g\in S^{m}$ we associate its Weyl quantization, namely
the operator $g^w(x,D_x)$, $D_x:=-\ir \partial_x$, defined by
\begin{equation}
\label{weyl}
G\psi(x)\equiv g^w(x,D_x)\psi(x):=\frac{1}{2\pi}\int_{\R^2}e^{\ir(x-y)\cdot
  \xi}g\left(\frac{x+y}{2};\xi\right) \psi(y)dyd\xi\ .
\end{equation} 
{\it We will often denote by a a capital letter the Weyl
  quantized of a symbol denoted with the corresponding lower case
  letter.} As an exception, we will denote by $W$ both the symbol of
the perturbation and the corresponding operator.

We use the symbol $\lambda(x,\xi)$  to {\it define,
  for $s\geq 0$ the
spaces $\cH^s=D([\lambda^w(x,-\ir \partial_x)]^{s(l+1)})$ (domain of the
$(s(l+1))^{th}$- power of the operator
operator 
$\lambda^w(x,-\ir \partial_x)$)
endowed by the graph norm}. For negative $s$, the space $\cH^s$ is the
dual of $\cH^{-s}$.

We will denote by $B(\cH^{s_1};\cH^{s_2})$ the space of bounded linear
operators from $\cH^{s_1}$ to $\cH^{s_2}$.

The potential $V$ defining 
$$
H_0:=H_\epsilon\big|_{\epsilon=0}\equiv \partial_{xx}+V
$$
is assumed to be a polynomial of
order $2l$, so that, in particular it belong to $S^{2l}$. 

We also assume that 
\begin{align}
\label{V3}
V'(x)\not=0\ ,\quad \forall x\not=0\ ,
\end{align}
and normalize the potential by assuming $V(0)=0$. 
The unperturbed Hamiltonian $H_0$ is the quantization of the
classical Hamiltonian system with Hamiltonian function
\begin{equation}
\label{h0c}
h_0(x,\xi):=\xi^2+V(x)\ .
\end{equation}
\begin{remark}
\label{periodic}
As a consequence of the assumptions above all the solutions of the
Hamiltonian system $h_0$ are periodic with a period $T(E)$ which
depends only on $E=h_0(x,\xi)$. 
\end{remark}

In the following we will denote by $\Phi^t_{h_0}$ the flow of the
Hamiltonian system \eqref{h0c}.

We denote by $\lambda_j^v$ the sequence of the eigenvalues of $H_0$ 
labeled in increasing order. It is well known that (see
e.g. \cite{HR82D})  
\begin{align}
\label{asyl}
\lambda_j\sim\frac{1}{c_l}j^d\ ,\quad j\to\infty
\end{align}
with $c_l>0$ and 
\begin{equation}
\label{d}
d=\frac{2l}{l+1}
\end{equation}
(in concrete examples one can compute also a complete asymptotic
expansion of the eigenvalues, see \cite{HR82D}). In the Harmonic case,
whithout lack of generality, we assume $V(x)=x^2$.

In what follows we will identify $L^2$ with $\ell^2$ by introducing
the basis of the eigenvector of $H_0$. Similarly we will identify
$\cH^{s}$ with the space  $\ell_s^2$ of the sequences
$\psi_j$ s.t. 
$$
\sum_{j\geq 1}j^{2s}\left|\psi_j\right|^2<\infty\ .
$$

In order to state the assumptions on the perturbation we need a few
notations. First we define the average with respect to the flow of
$h_0$:
\begin{equation}
\label{mediaW}
\langle W\rangle(x,\xi,\omega t):=\frac{1}{T(E)}\int_0^{T(E)}W\left(
\Phi^\tau_{h_0}(x,\xi),\omega t \right)d\tau\ .
\end{equation}

Concerning the perturbation, we assume that $W\in
\cS^{\beta}$ and we define
\begin{equation}
\label{betatilde}
\tilde\beta:=\left\{
\begin{matrix}
2\beta-2l & {\rm if} &\langle W\rangle\equiv
0
\\
\beta &\null &{\rm otherwise}
\end{matrix}
\right.\ .
\end{equation}

The main result of the paper is the following theorem.
\begin{theorem}
\label{m.1}
Assume $\tilde\beta<l+1$, then there exist $\epsilon_*>0$, $C_\lambda$
and $\forall \left|\epsilon\right|<\epsilon_*$ a closed set
$\Omega(\epsilon)\subset\Omega$ and, $\forall
\omega\in\Omega(\epsilon)$ there exists a unitary (in $L^2$) time
quasiperiodic operator $\Phi_\omega(\omega t)$ s.t. the function $\varphi$
defined by $\Phi_\omega(\omega t)\varphi:=\psi $ satisfy the equation
\begin{equation}
\label{rido}
\ir \dot \varphi= H_{\infty}\varphi\ ,
\end{equation}  
with $H_{\infty}=\diag (\lambda_j^\infty)$ and 
\begin{equation}
\label{dla}
\left|\lambda_j^\infty-\lambda_j^v\right|\leq C_{\lambda}\epsilon
j^{\frac{\tilde \beta}{l+1}}\ .
\end{equation}
Furthermore one has
\begin{itemize}
\item[1.] $\displaystyle{\lim_{\epsilon\to0}}\left|\Omega-\Omega(\epsilon)\right|=0$;
\item[2.] $\forall s$ $\exists \epsilon_s$ s.t., if
  $\left|\epsilon\right|<\epsilon_s$ then $\Phi_\omega(\omega t)\in
  B(\cH^s;\cH^{s})$;
\item[3.] $\forall r>0$ $\exists \epsilon_{s,r}>0$ s.t. if
  $\left|\epsilon\right|<\epsilon_{s,r}$ then $\exists s_r$ s.t. the map
  $\phi\mapsto \Phi_\omega(\phi)$ is of class
  $C^r(\T^n;B(\cH^{s+s_r};\cH^{s}))$;
\item[4.] $\exists a>0$, s.t.
  $\norma{\Phi_\omega(\phi)-\uno}_{B(\cH^{s+\beta};\cH^{s})
}\leq C_s\epsilon^a$.
\end{itemize}
\end{theorem}
\begin{remark}
\label{r.main1}
Under the assumptions of the Theorem \ref{m.1}, the perturbation $W$ is
an unbounded operator; it is for this reason that $\Phi_\omega$ is
close to identity only as an operator decreasing smoothness.
\end{remark}
\begin{remark}
\label{r.main2}
With our technique we are not able to show that the sequence
$\epsilon_r$ does not go to $0$ as $r\to\infty$, thus we cannot
guarantee that $\Phi_\omega$ is actually a $C^\infty$ function of the angles.  
\end{remark}
\begin{remark}
\label{r.main.3}
The dependence of $\Phi_\omega$ on $\omega$ is Whitney smooth; however, for the
sake of simplicity we did not work out a precise statement. 
\end{remark}

A consequence of the above theorem is that in the considered range of
parameters all the Sobolev norms, i.e. the $\cH^s$ norms of the
solutions are bounded forever and the spectrum of the Floquet operator
is pure point.

A couple of examples is useful in order to clarify the range of
applicability of the result.

\begin{example}
\label{ex.uno}(Duffing oscillators)
$l=2$. The assumptions of Theorem \ref{m.1} become $\beta<3$ if
$\langle W\rangle\not =0$ otherwise $\beta< 7/2$. An example in which
the assumption are fulfilled is a singular version of the 
Duffing oscillator:
\begin{equation}
\label{duffing}
-\partial_{xx}+x^4+\epsilon x^\beta f(\omega t)\ ,\ \beta=1,2,3
\end{equation}
where $f$ is an arbitrary $C^\infty$ function. (In this case one
has that, for symmetry reasons the average of $x^3$ is zero.)
At the end of the section we will show that the method of the present paper
can be extended to deal also with the case $\beta=4$. The best
previous result, due to \cite{LY10}, only allowed to have $\beta=1$.

One can also add a magnetic type term of the form
\begin{equation}
\label{mag.1}
-(a_0(\omega t)+a_1(\omega t)x)\ir \partial_x \ .
\end{equation}
More general perturbations of the form of a pseudodifferential
operator with symbol $W(x,\xi,\omega t)$ with $W$ of class
$S^{\beta}$, $\beta<3$. 
are allowed.
\end{example}

\begin{example}
\label{harmonic}
(Harmonic oscillator) $l=1$. In this
case the Theorem \ref{m.1} applies when $\beta<2$. Thus, for example we
can deal with the case
\begin{equation}
\label{har.1}
-\partial_{xx}+x^2+\epsilon x a_1(\omega t)-\ir a_2(\omega
t)\epsilon\partial_x \ .
\end{equation}
The more general case of a perturbation quadratic in $x$ and
$\xi$ will be covered at the end of the section.

Perturbations of the kind of those considered by Delort \cite{D14}
belong to the class of symbols dealt with in Theorem \ref{m.1}. The
same is true for the main term of the perturbation in \cite{GY00}.
\end{example}

\begin{remark}
\label{nonpoly}
If $W\in S^\beta$ is independent of $\xi$, namely $W=W(x)$, then it
must be a polynomial. Indeed, if $k>\beta$, then $|\partial_x^kW(x)|$
must tend to zero as $\xi\to\infty$, and thus it must be identically
zero. 

As anticipated in the introduction, the extension of Theorem \ref{m.1}
to more general perturbations including the cases of the form
\eqref{nonp} with $a_0$, $a_1$ non polynomial smooth functions will be
obtained in paper II.
\end{remark}

In order to give the extension to $\beta=2l$ (and also for future use)
it is useful to give the definition of quasihomogeneuos symbols.

\begin{definition}
\label{quasi.hom.d}
We will say that a symbol $f$ is quasihomogeneous of degree $m$
if
\begin{equation}
\label{q.hom}
f(\rho x,\rho^l\xi)=\rho^m f(x,\xi)\ ,\quad \forall \rho>0\ .
\end{equation}
\end{definition}

The most general time dependent quasihomogeneous polynomial of degree
$2l$ is given by
\begin{equation}
\label{q.h.2l}
W_{2l}(x,\xi,\omega t)=a_1(\omega t)\xi^2+a_2(\omega t)x\xi+a_3(\omega
t) x^{2l}\ .
\end{equation}

\begin{theorem}
\label{quasi.l}
Consider the Schr\"odinger
equation with Hamiltonian 
\begin{equation}
\label{h.l}
-\partial_{xx}+x^{2l}+\epsilon W_{2l}(x, D_x,\omega t)+\epsilon
W(x,D_x,\omega t)\ ,
\end{equation}
with $l\geq 1$ a positive integer and $a_j\in C^\infty(\T^n)$, and $W\in
\cS^\beta$ fulfilling the assumptions of Theorem \ref{m.1}, then the
same conclusions of Theorem \ref{m.1} hold (in eq. \eqref{dla} one has
to put $\tilde \beta=l+1$)
\end{theorem}

\section{Transformations of linear time dependent equations}\label{algo}

In the following we will use in some different contexts
transformations of the form $\psi=e^{-\ir\epsilon X(\omega t)}\vf$, with $X$
a family of self adjoint operators that in some sense depend smoothly
on time. So, to start with, we study in a purely formal way how the
Schr\"odinger equation is changed by such  transformations. In the
subsequent sections we will make all notions precise.

\begin{definition}
\label{lie}
Let $X$ be a selfadjoint operator; we will say that 
\begin{equation}
\label{qlie}
(Lie_{\epsilon X}F):=\meX F\eX  
\end{equation}
is the quantum Lie transform of $F$ generated by $\epsilon X$.
\end{definition}
Remark that the quantum Lie transform fulfills the equation
$$
\frac{d}{d\epsilon}Lie_{\epsilon X}F=-\ir\left[Lie_{\epsilon X}F;X\right]= \meX \ir
     [X;F]\eX\ ,
$$
from which one immediately gets (formally!)
\begin{align}
\label{qlie1}
 Lie_{\epsilon X}F =\sum_{k\geq 0}\frac{1}{k!}\epsilon^k F_k\ ,
\\
\label{serqlie}
F_0=F\ ;\quad F_k:=-\ir [F_{k-1};X]\ .
\end{align}
Remark also that one has
\begin{equation}
\label{r.lie}
\frac{d^k}{d\epsilon^k}Lie_{\epsilon X}F=\meX F_k\eX\ .
\end{equation}

In the following we will meet situations where the above series are
either convergent or asymptotic. 

We will use the same terminology also when $X$ depends on time and/or
on 
$\omega$ (which in this case play the role of parameters).

\begin{lemma}
\label{T.1}
Let $F$ be a selfadjoint operator, and let
$X(t)$ be a family of selfadjoint operators. Assume that
$\psi(t)$ fulfills the equation
\begin{equation}
\label{1}
\im\dot \psi=F\psi\ ,
\end{equation}  
then $\vf$ defined by 
\begin{equation}
\label{2}
\vf=e^{\im\epsilon X(t)}\psi\ ,
\end{equation}
fulfills the equation 
\begin{equation}
\label{3}
\im\dot \vf =F_\epsilon(t) \vf 
\end{equation}
with 
\begin{align}
\label{4.1.1}
F_\epsilon&:=Lie_{\epsilon X} F-Y_X\ ,
\\
\label{yx}
&Y_X:=\int_0^\epsilon (Lie_{(\epsilon-\epsilon_1) X}\dot X)d\epsilon_1\ .
\end{align}
\end{lemma}
\proof One has 
$$
\ir\frac{d\vf}{dt}=\ir\frac{d \meX}{dt}\psi+\meX F\eX\vf
= \left(\ir\frac{d \meX}{dt}\eX+\meX F\eX
\right)\vf\ .
$$ 
So, the second term in the bracket is already $Lie_{\epsilon X}F$. Define 
$$
\tilde Y_X:=\frac{d \meX}{dt}\eX\ ,
$$
and compute 
\begin{align*}
\frac{d \tilde Y_X}{d\epsilon}=\frac{d}{dt}\left(\im X\meX \right)\eX-\im 
\frac{d \meX}{dt}\eX X
\\
=\im \dot X+\im X\frac{d\meX}{dt} \eX-\im \frac{d\meX}{dt}\eX
X
\\
=
\im\dot X -\im\left[\tilde Y_X,X \right]
\end{align*}
It follows that $\tilde Y_X$ solves the Cauchy problem
\begin{align*}
\frac{d\tilde Y_X}{d\epsilon}=\ir \dot X-\ir [\tilde Y_X;X]\ ,\quad
\tilde Y_X(0)=0\ ,
\end{align*}
whose solution is easily computed by Duhamel formula getting
\eqref{yx}.
\qed

\begin{definition}
\label{t.tras}
Given $X$, we will say that 
\begin{equation}
\label{trasf.X.1}
T_{\epsilon X}F:= Lie_{\epsilon X}F-Y_{X} 
\end{equation}
is the transformation of $F$ through $\epsilon X$. Remark that 
$$
T_{\epsilon X}(F+G)=T_{\epsilon X}F+ Lie_{\epsilon X}G\ .
$$ 
\end{definition}

\begin{remark}
\label{first}
In the following we will be interested in expansions either in
$\epsilon$ or in operators which are more and more regularizing; in
this second case, as usual, the key property that we use is that the
commutator of two operators is more regularizing than the product of
the original operators. Thus, up to higher order corrections, either
in $\epsilon$ or in smoothness, we will have that if $F$ has the
structure $F=H_0+\epsilon P$ with $P$ more smoothing (or ``less
unbounded'') than $H_0$, then, up to higher order corrections, one has 
\begin{equation}
\label{primo}
T_{\epsilon X}F=H_0+\epsilon P-\ir \epsilon[ H_0;X]- \epsilon\dot X+... 
\end{equation}
\end{remark}

\section{Smoothing the perturbation}\label{egor}

\subsection{Some symbolic calculus}\label{symbol}

First we recall that, from the Calderon Vaillencourt Theorem, the
following lemma holds.
\begin{lemma}
\label{caderon}
Let $f\in S^{m}$, then one has
\begin{equation}
\label{CV}
f^w(x,D_x)\in B(\cH^{s_1+s};\cH^{s})\ ,\quad \forall s\ ,\quad \forall
s_1\geq m\ .
\end{equation}
\end{lemma}
We emphasize that the result holds also for negative values of the
indexes $m,s_1$. 

Given a symbol $g\in S^{m}$ we will write
\begin{equation}
\label{asym}
g\sim\sum_{ j\geq 0}g_j\ ,\quad g_j\in
S^{m_j}\ ,\quad m_j\leq
m_{j-1}\ ,
\end{equation}
if $\forall \kappa$ there exist $N$ and $r_N\in
S^{-\kappa}$ s.t.
$$
g=\sum_{j=0}^{N}g_j+r_N\ .
$$

The following result is standard 
\begin{lemma}
\label{l.M}
Given a couple of symbols $a\in S^{m}$ and $b\in S^{m'}$, denote by
$a^w(x,D_x)$ and $b^w(x,D_x)$ the corresponding Weyl
operators, then there exists a symbol $c$, denoted by $c=a\sharp b$ such
that
$$
(a\sharp b)^{w}(x,D_x)=a^w(x,D_x)b^w(x,D_x)\ ,
$$ 
furthermore one has
\begin{equation}
\label{sharp}
(a\sharp b)\sim \sum_{j} c_j
\end{equation}
with 
$$
c_j=\sum_{k_1+k_2=j}\frac{1}{k_1!k_2!}\left(\frac{1}{2}\right)^{k_1}
\left(-\frac{1}{2}\right)^{k_2}  (\partial^{k_1}_\xi D^{k_2}_xa)
(\partial^{k_2}_\xi D^{k_1}_xb)\in S^{m+m'-(l+1)j} \ .
$$
In particular, denoting\footnote{Sometimes $\left\{.;.\right\}^q$ is
  called the Moyal Bracket} 
$$\left\{a;b\right\}^q:=-\ir(a\sharp b-b\sharp
a)\ ,$$ we have
\begin{equation}
\label{moy}
\left\{a;b\right\}^q=\left\{a;b\right\}+S^{m+m'-3(l+1)}\ ,
\end{equation}
where 
$$
\left\{a;b\right\}:=-\partial_\xi a\partial_xb+\partial_\xi
b\partial_xa\in S^{m+m'-(l+1)}\ ,
$$
is the Poisson Bracket between $a$ and $b$, while \eqref{moy} means
that $\left\{a;b\right\}^q=\left\{a;b\right\}+$some quantity belonging
to  $\in S^{m+m'-3(l+1)}$. {\it Similar notations will be systematically
used in the following.}
\end{lemma}

Sometimes we will deal with symbols having finite differentiability. We
will denote by $S^{m}_N$ the space of symbols which are only $N$ times
differentiable and fulfill the inequality \eqref{sm} only for
$k_1+k_2\leq N$. This is a Banach space with the norm
\begin{equation}
\label{Sn}
\norm{g}_{S^{m}_N}:=\sum_{k_1+k_2\leq
  N}\sup_{(x,\xi)\in\R^2}\frac{\left|\partial^{k_1}_x
  \partial_{\xi}^{k_2}g(x,\xi)
  \right|}{[\lambda(x,\xi)]^{m-lk_1-k_2}}\ .
\end{equation}

We remark that for the space $\cS^{m}$ a family of seminorms is given by
the standard norms of $C^M(\T^n;S^{m}_N)$ as $M$ and $N$ vary.

Finally we will deal with Whitney smooth functions of the
frequencies\footnote{This will be needed only for the proof of Lemma
  \ref{diaga.300}. For the rest of KAM theory Lipschitz dependence on
  the frequencies is enough.}. To this end we recall (following
\cite{Stein}) the definition of smooth function on a closed set
$\widetilde \Omega\subset \Omega$. Fix an integer $k$ and a $\rho$
fulfilling $k<\rho\leq k+1$; let $\cB$ be a Banach space, and
$f:\widetilde \Omega\to\cB$ a map. The map $f$ is said to be of class
$Lip_\rho(\widetilde\Omega;\cB)$, if there exist maps $f^{(j)}$,
$0\leq|j|\leq k$ defined on $\widetilde\Omega$, such that $f^{(0)}=f$
and so that, if
\begin{equation}
\label{why.1}
f^{(j)}(\omega)=\sum_{|j+l|\leq k}\frac{f^{(j+l)}(\nu)}{l!}(\omega-\nu)^l
+R_j(\omega,\nu)\ , 
\end{equation}
then
\begin{equation}
\label{why.2}
\left\|f^{(j)}(\omega)\right\|\leq M\ ; \quad \left\|
R_j(\omega,\nu)\right\|\leq M\left|\omega-\nu\right|^{\rho-k}\ ;\quad
\forall \omega,\nu\in\widetilde \Omega\ ,\quad \left|j\right|\leq k\ .
\end{equation}
Here we used a standard vector notation: $j=(j_1,...,j_n)$ and
$\omega^j\equiv \omega_1^{j_1}...\omega_{n}^{j_n}$. The minimum of the
constants $M$ for which \eqref{why.2} holds is a norm on the space  $Lip_\rho(\widetilde\Omega;\cB)$

\begin{definition}
\label{why.1.111}
We will say that a function $f:\widetilde\Omega\to\cS^{m}$ is of
class $Lip_\rho^{m}(\widetilde\Omega)$ if for all $N_1,N_2$ it is of
class $Lip_\rho(\widetilde\Omega;C^{N_1}(\T^n;S^{m}_{N_2}))$.
\end{definition}

\begin{definition}
\label{pseud}
An operator $F$ will be said to be a pseudodifferential operator of
class $O^{m}$ if there exists a sequence $f_j\in
S^{m_j}$ with $m_j\leq
m_{j-1}$ and, for any $\kappa$ there exist $N$ and
an operator $R_N\in B(\cH^{s-\kappa};\cH^{s})$, $\forall s$ such that
\begin{equation}
\label{expa}
F=\sum_{j\geq0}^Nf_j^w+R_N\ .
\end{equation}
In this case we will write $f\sim\sum_{j\geq 0} f_j$ and $f$ will be
said to be the symbol of $F$.
\end{definition}

Concerning maps we will use the following definition

\begin{definition}
\label{pseudomap}
A map $\T^n\ni\phi\mapsto F(\phi)\in O^{m}$, will be said to be
of class $\cO^{m}$ if the functions of the sequence $f_j$
also depend smoothly on $\phi$, namely
$f_j\in\cS^{m_j}$ and the operator valued map
$\phi\mapsto R_N(\phi)$ has the property that for any $K\geq1 $ there
exists $a_K\geq 0$ s.t. for any $N$ one has
\begin{equation}
\label{propcN}
R_N(.)\in C^K(\T^n;B(\cH^{s-\kappa+a_K};\cH^{s}))\ , \forall s \ .
\end{equation} 
\end{definition}

\begin{definition}
\label{why.4}
A map $\tilde\Omega\ni\omega\mapsto F\in\cO^{m}$ will be said to be of
class $\cL ip_\rho^{m}(\widetilde \Omega)$ if the functions $f_j\in
Lip_\rho^{m_j}(\widetilde\Omega)$ and if the map $\omega\mapsto R_N$
has the property that there exists $b\geq0$ s.t.
\begin{equation}
\label{lip}
R_N\in Lip_\rho (\tilde \Omega;
C^K(\T^n;B(\cH^{s-\kappa+a_K+b};\cH^{s})))\ , \forall s \ .
\end{equation}
\end{definition}

We want now to study the quantum Lie transform generated by a symbol
$\chi\in S^{m}$. First, applying Proposition A.2 of \cite{MR16} we
have the following Lemma  
\begin{lemma}
\label{MR}
Let $\chi\in S^{m}$ with $m\leq l+1$, then $X:=\chi^w(x,D_x)$ is
selfadjoint and $e^{-\ir \epsilon X}$ leaves invariant all the spaces
$\cH^s$.
\end{lemma}
\proof According to \cite{MR16}, the thesis holds if there exists a positive
selfadjoint operator $K$ such that both the operators $XK^{-1}$ and
$[X,K]K^{-1}$ are bounded. To this end we take $K$ to be the Weyl
operator of the symbol $\lambda^{m}$. From symbolic calculus it
follows that $XK^{-1}\in O^{0}$ and
$[X,K]K^{-1}\in O^{2m-(l+1)-m}$. Thus they are bounded under the
assumption of the Lemma.\qed

One can rewrite formulae \eqref{qlie1} and \eqref{serqlie} in terms
of symbols. Thus, if $f$ and $\chi$ are symbols and $\chi$ fulfills the
assumptions of lemma \ref{MR} one can define
\begin{equation}
\label{liqs}
f_0^q:=f\ ,\quad f_k^q:=\left\{ f_{k-1}^q;\chi\right\}^q\ ,
\end{equation}
and one can expect that the symbol of $Lie_{\epsilon X}F$ is $\sum_{k\geq
  0}\epsilon^kf_k^q/k!$. This is ensured by the following lemma: 

\begin{lemma}
\label{egorov}
Let $\chi\in S^{m}$ and $f\in S^{m'}$ be symbols, assume that $m<l+1$,
then $Lie_{\epsilon X}F\in O^{m}$, and furthermore its symbol, denoted by
$lie_{\epsilon\chi} f$ fulfills
\begin{equation}
\label{liqser}
lie_{\epsilon\chi} f\sim \sum_{k\geq
  0}\frac{\epsilon^kf_k^q}{k!} \ .
\end{equation}
\end{lemma}
\proof First remark that, by induction, one has $f^q_k\in
S^{m'+k(m-l-1)}$. From \eqref{r.lie} and the formula
of the remainder of the Taylor expansion one also has
$$
Lie_{\epsilon X}F(\epsilon)=\sum_{k=0}^{N}\frac{F_k}{k!} 
\epsilon^k+\frac{\epsilon^{N+1}}{N!}
\int_0^1(1+u)^J e^{-\ir u\epsilon X}F_{N+1}e^{\ir
  u\epsilon X}du\ ,
$$
so that, by defining $R_N$ to be the integral term of the previous
formula, we have $R_N\in B(\cH^{s-\kappa},\cH^s)$ with
$\kappa=m'+(m-l-1)N$, which diverges as $N\to\infty$
and thus shows that the expansion \eqref{liqser} is asymptotic in the
sense of definition \ref{pseud}.\qed 

In the following, by abuse of language, we will call
$lie_{\epsilon\chi} f$ the quantum Lie transform of $f$ through $\chi$.

\begin{remark}
\label{eg.1}
Denote by $\Phi^\epsilon_\chi$ the flow of the Hamilton equations of
$\chi$, then one has 
\begin{align}
\label{liec}
f\circ \Phi^\epsilon_\chi \sim \sum_{k\geq
  0}\frac{\epsilon^kf_k}{k!} 
\\
\label{liec.1}
f_0:=f\ ,\quad f_k:=\left\{ f_{k-1};\chi\right\} \ ,
\end{align}
thus one has
\begin{align}
\label{eg.3}
lie_{\epsilon\chi} f= f\circ
\Phi^\epsilon_\chi+S^{m+m'-3(l+1)}
\\
\label{eg.3.1}
=f+\left\{f;\chi\right\} +S^{m'+2(m-l-1)}\ .
\end{align}

\end{remark}
In the following we will need also a result valid in the limit case
$\chi\in S^{l+1}$. This is covered by the following Lemma, which is a
variant of Theorem 7.1 of \cite{HR82}.

\begin{theorem}
\label{7.1}
Let $\chi\in S^{l+1}$ and
let $f\in S^{m'}$; assume that $f\circ \Phi_\chi^\epsilon\in
S^{m'}$, then equation \eqref{eg.3} holds. 
\end{theorem}
The proof is obtained exactly as in \cite{HR82} and is omitted.

\begin{remark}
\label{yx.1}
Let $\chi\in S^m$ then the operator $Y_X$ defined by
eq. \eqref{yx} is a pseudodifferential operator with symbol
\begin{equation}
\label{yx.2}
y_x:=\int_0^\epsilon(lie_{(\epsilon-\epsilon_1)\chi}\dot\chi)d\epsilon_1=
\epsilon\dot\chi +\epsilon S^{ 2m-(l+1)}\ . 
\end{equation}
\end{remark}

\subsection{Symbol of the transformed Hamiltonian and formal
  description of the smoothing algorithm.}\label{smoa}

The idea is to use the quantization of a time dependent symbol $\chi(\omega
t)$ in order to transform the original Hamiltonian 
\begin{equation}
\label{ts.1}
h:=h_0+\epsilon W
\end{equation}
into a new one with a more regular perturbation. 

According to eq. \eqref{primo}, written at the
level of symbols, one has that the tranformed
Hamiltonian has a symbol which, at highest order, is given by
\begin{equation}
\label{primo.1}
h_\epsilon=h_0+\epsilon W+\epsilon\left\{h_0;\chi\right\}^q-
\epsilon\dot\chi +...=h_0+\epsilon W+\epsilon\left\{h_0;\chi\right\}-
\epsilon\dot\chi +...
\end{equation}
So, in order to increase the order of the perturbation one has to
choose $\chi$ in such a way to eliminate the terms of order $\epsilon$
or to transform them into smoother objects. To explain the
procedure one has to distinguish between the case $l>1$ and the case
$l=1$. 

Consider first $l>1$. In this case it turns out that $\dot \chi$ is
more regular then $\left\{h_0;\chi\right\}$ (see Lemma \ref{lemchi}),
so in that case one determines $\chi$ by solving the homological
equation
\begin{equation}
\label{12}
p+\left\{h_0;\chi\right\}=\langle p\rangle\ .
\end{equation}
with $p=W$ (this will be done in Lemma \ref{lemchi}).
Using such a $\chi$ to transform the Hamiltonian, one gets a new
Hamiltonian with a symbol which is a perturbation of
\begin{equation}
\label{pr.1}
h_0+\epsilon\langle W\rangle(h_0,\omega t)\ ,
\end{equation}
Remark that $\langle W\rangle$ is a function of the phase variables
only through $h_0$ (since it Poisson commutes with it), but it is also
time dependent.

So, the second step consists in looking for a second
generating function $\chi_1=\chi_1(h_0,\omega t)$ in order to
eliminate the time dependence from $\langle W\rangle$ (at the main
order). Taking into account that in such a case
$\left\{h_0;\chi_1\right\}\equiv 0$, the main term of the Hamiltonian
transformed through such a $\chi_1$ is simply given by
$$
h_0-\epsilon \dot \chi_1+\epsilon\langle W\rangle(h_0,\omega t)+...
$$ and this leads to the second kind of homological equation that we
need to solve, namely
\begin{equation}
\label{12.p}
-\omega\cdot\frac{\partial \chi_1}{\partial\phi}=p-\bar p\ ,
\end{equation}
where $p\equiv\langle W\rangle $, while $\bar p$ is defined by
\begin{equation}
\label{media.2}
\bar p(x,\xi):=\frac{1}{(2\pi)^{n}}\int_{\T^n}p(x,\xi,\phi)d\phi\ .
\end{equation}

Using such a $\chi_1$ one transforms the Hamiltonian into a
perturbation of 
\begin{equation}
\label{sp.3}
h_0+\epsilon\overline{\langle W\rangle}(h_0)
\end{equation}
which is a function of $h_0$ only. Thus the idea is to repeat the
procedure with $h_0$ replaced by the function \eqref{sp.3}. As a
consequence, at the subsequent steps, we will have to solve homological
equations of the form of \eqref{12} with $h_0$ replaced by a function
of $h_0$, and this will lead to the homological equation 
\begin{equation}
\label{12.3}
p+\left\{h_1;\chi\right\}=\langle p\rangle\ ,
\end{equation}
with
\begin{equation}
\label{h1}
h_1:=h_0+\epsilon f(h_0)\ ,
\end{equation}
which will be solved thanks to the Remark \ref{lemchi.1}.  Then one
can proceed interatively until the perturbation is reduced to a
smoothing operator of arbitrary order. Actually the procedure we use
is slightly modified in order to be able to deal with a singularity
related to the singularity of the action variables at the origin and
in order to get a better result when the average of $W$ vanishes (see
the proof of Theorem \ref{smoothing}).
 
In the case $l=1$ the situation is different since in this case
$\dot\chi$ and $\left\{ h_0;\chi\right\}$ belong to the same
smoothness class. So in this case we consider again equation
\eqref{primo.1}. In order to reduce all the terms of order $\epsilon$
one has to solve the following homological equation
\begin{equation}
\label{12.1}
\left\{h_0,\chi\right\}-\dot\chi+W=\overline{\langle W\rangle}\ ,
\end{equation}
and in this case the original Hamiltonian is directly transformed into
a new one of the form
\begin{equation}
\label{sp.4}
h_0+\epsilon\overline{\langle W\rangle}(h_0)+\epsilon p(x,\xi)
\end{equation}
with $p$ which is the symbol of a more smoothing operator. We remark
that equation \eqref{12.1} can only be solved in the case where the
period of the orbits of $h_0$ does not depend on the energy, and
therefore only in the case where $V$ is exactly quadratic.  Now there
is a difficulty: one cannot include $\overline{\langle W\rangle}$ in
the main part of the Hamiltonian in order to iterate since this would
eliminate the above property. However it turns out that this is not
needed, since in the Harmonic case one has
$\left\{h_0,f\right\}^q=\left\{h_0,f\right\}$. This allows to proceed
as in classical normal form theory and to conjugate the Hamiltonian to
a very smoothing symbol.

\subsection{Solution of the homological equations}\label{l.s.1}

From now on we will use the notation 
\begin{equation}
\label{sleq}
a\sleq b
\end{equation}
to mean ``there exists a constant $C$ independent of all the relevant
quantities, such that $a\leq Cb$''.

In the following we will meet functions which depend on the phase
space variables only through $h_0$, namely functions $p$ such that
there exist a $\tilde p$ with the property that
$$
p(x,\xi)=\tilde p(h_0(x,\xi))\ .
$$
For such functions we introduce a new class of symbols. 

\begin{definition}
\label{d.sm}
A function $\tilde p\in\C^{\infty}$ will be said to be of class
$\widetilde S^m$
if one has
\begin{equation}
\label{sm.1}
\left|\frac{\partial^k\tilde p}{\partial E^k}(E)\right|\sleq \langle
E^{\frac{m}{2l}-k}\rangle \ .
\end{equation} 
\end{definition}

We will also need to use functions from $\T^n$ to $\widetilde S^m$
which may also depend in a (Whitney) smooth way on the
frequencies. For these classes we will use the same notation we
already introduced, simply we will put a tilde on the letter denoting
the corresponding class. Furthermore, by abuse of notation, we will
say that $p\in\widetilde S^m$ if there exists $\tilde p\in \widetilde
S^m$ s.t. $p(x,\xi)=\tilde p(h_0(x,\xi))$.

As one can see in the case of a homogeneous potential $V(x)=x^{2l}$,
$l>1$, the period as a function of the energy has a singularity at
zero. In order to avoid this problem, before starting the procedure,
it is useful to modify the perturbation making a cutoff close to the
origin.

Let $\eta$ be a $C^\infty$ function such that
\begin{equation}
\label{cutoff}
\eta(E)=\left\{   
\begin{matrix}
1 & {\rm if} & |E|>2
\\
0 & {\rm if} & |E|<1
\end{matrix}
\right.
\end{equation}
and split 
\begin{equation}
\label{Wsplit}
W=W_0+W_\infty\ ,\quad W_\infty(x,\xi)=W(x,\xi)(1-\eta(h_0(x,\xi)))\ ,\quad
W_0(x,\xi)=W(x,\xi)\eta(h_0(x,\xi))\ , 
\end{equation}
then $W_\infty\in \cS^{-\kappa}$ for any $\kappa$, and $W_1\in
S^\beta$ is the actual perturbation that has to be regularized.

\begin{remark}
\label{locality}
All the smoothing procedure is based on the solution of the
homological equation and computation of Moyal brackets, which (up to
operators which are smoothing of all orders) are operations preserving the
property of symbols of being zero in the region $E<1$. 
\end{remark}

\begin{lemma}
\label{period}
Consider the period $T(E)$. Then the function $\bar T(E):=\eta(E)T(E)$
is a symbol and one has $\bar T\in \widetilde
S^{1-l}$. 
\end{lemma}
\proof Consider the function
$$
A(E):=\eta(2E)\int_{\left\{ h_0(x,\xi)\leq E \right\}} dxd\xi 
$$ According to Lemma (1-3) of \cite{HR82D} this is a symbol of class
$\tilde S^{l+1}$. But this function, when $E>1$ is the classical
action of the Hamiltonian system $h_0$. Thus, in the
region $E>1$, one has $T(E)=2\pi\partial A/\partial E$. Now, $\bar T$
coincides with this function in the considered region and is regular
and bounded in the other region, and thus the thesis follows.\qed

\begin{remark}
\label{per.2pi}
The function $A(E)$ is particularly important since the period of the
orbits of the Hamiltonian system $A(h_0)$ is $2\pi$ whenever
$E>1$. Furthermore, exploiting the fact that in the region $E>1$
$A(h_0)$ admits an expansion in quasihomogeneous polynomials (see the
Appendix of \cite{HR82}), one can see that given a symbol $f\in S^m$
then $f\circ\phi^t_{A(h_0)}\in S^m$.
\end{remark}

\begin{lemma}
\label{lememd}
Let $p\in S^{m}$ be a symbol supported in the region $h_0(x,\xi)>1$,
then $\langle p\rangle\in\widetilde S^m$.
\end{lemma}
\proof Consider the function $A(h_0)$. It is easy to see that, in the region
$E>1$, 
\begin{equation}
\label{ah0}
\Phi^t_{h_0}=\Phi^{\frac{t2\pi}{\bar T}}_{A(h_0)}\ .
\end{equation}
therefore one has 
\begin{align*}
\frac{1}{T(h_0)}\int_0^{T(h_0)}p\circ
\Phi^t_{h_0}dt =\frac{1}{\bar T(h_0)}\int_0^{\bar T(h_0)}p\circ
\Phi^t_{h_0}dt=\frac{1}{\bar T(h_0)}\int_0^{\bar T(H_0)}p\circ
\Phi^{2\pi t/\bar
  T}_{A(h_0)}dt
\\
=\frac{1}{2\pi}\int_0^{2\pi}p\circ\Phi^t_{A(h_0)}dt  
\ ,
\end{align*} 
but $p\circ
\Phi^t_{A(h_0)}\in S^m$, so that the result immediately follows.\qed

Concerning the solution of the homological equation \eqref{12} we have
the following Lemma.

\begin{lemma}
\label{lemchi}
Let $p\in S^{m}$ be a symbol which vanishes in the region $h_0<1$,
then the homological equation \eqref{12} has a solution $\chi$ which
is a symbol of class $\chi\in S^{m-l+1}$.
\end{lemma}
\proof First, following Lemma 5.3 of \cite{BG93}, we have that $\chi$
is given by the formula
\begin{equation}
\label{coh.2}
\chi=\frac{1}{T(E)}\int_0^{T(E)}t\check p\circ \Phi^t_{\hz} dt\ ,
\end{equation}
with $\check p:=p-\langle p\rangle$. To see this, fix a value of $E$ and
compute
\begin{align*}
\left\{\chi;{h_0}\right\}(\zeta)&=\left.\frac d{dt}\right|_{t=0} 
\chi\left(\Phi^t_{h_0}(\zeta)\right)= \left.\frac d{dt}\right|_{t=0} \frac1T 
\int_0^T\check p\left(\Phi^{t+s}_{h_0}(\zeta)\right)sds
=\frac1T\int_0^Ts\left.
\frac 
d{ds}\check p\left(\Phi^{t+s}_{h_0}(\zeta)\right)ds\right|_{t=0}
\\
&=\frac1T\check p\left(\Phi_{h_0} 
^{s}(\zeta)\right)s\big|_{0}^T-\frac1T\int_0^T\check p\left(\Phi^{s}_{h_0}(\zeta) 
\right)ds=\check p(\zeta)\ ,
\end{align*}
where $\zeta=(x,\xi)$. Now, exploiting again \eqref{ah0}, one has
\begin{equation*}
\chi=\frac{1}{T(E)}\int_0^{T(E)} t\check p\circ \Phi^t_{\hz}
dt=\frac{1 }{\bar T(E)}\int_0^{\bar
  T(E)}t\check p\circ\Phi^{\frac{t 2\pi}{\bar
    T(E)}}_{A(h_0)}dt 
=\frac{\bar
  T(E)}{4\pi^2}\int_0^{2\pi}t\check  p\circ\Phi^{t}_{A(h_0)}dt \ ,
\end{equation*}
from which, exploiting Lemma \ref{period}, one immediately gets
the result. \qed
\begin{remark}
\label{lemchi.1}
From the above proof one  gets that the above technique
also allows to solve the homological equation \eqref{12.3} and to show
that the solution also belongs to $S^{m-l+1}$. 
\end{remark}
\begin{remark}
\label{funzioni}
In the above lemmas $p$ can also depend on the angles $\phi$ and on
the frequencies $\omega$, but they only play the role of parameters,
so in that case the result is still valid substituting the classes
$\cS$ or $Lip_\rho$ to the classes $S$ with the same index.
\end{remark}

We come now to equation \eqref{12.p}. 

First, fix $\tau>n-1$ and denote 
\begin{equation}
\label{diof}
\Omega_{0\gamma}:=\left\{\omega\in\Omega\ :\ \left|k\cdot\omega\right|\geq\gamma
|k|^{-\tau}\right\}\ ,
\end{equation}
then it is well known that 
\begin{equation}
\label{sti}
\left|\Omega-\Omega_{0\gamma}\right|\sleq \gamma\ .
\end{equation}

\begin{lemma}
\label{coh2}
Let $p\in \widetilde{Lip}_\rho^{m}(\Omega_{0\gamma})$, then there
exists a solution $\chi\in \widetilde{Lip}_\rho^{m}(\Omega_{0\gamma})
$ of \eqref{12.p}. Furthermore, in this case $\bar p\in
\widetilde{Lip}_\rho^{m}(\Omega_{0\gamma}) $.
\end{lemma}
\proof We proceed as usual expanding $p$ in Fourier series. First we
consider the case where $p$ does not depend explicitly on $\omega$.
Define
\begin{equation}
\label{four}
p_k(E):=\frac{1}{(2\pi)^n}\int_{\T^n}p(E,\phi)e^{-\im k\phi}d\phi\ ,
\end{equation}
and remark that, since $\forall MN$ the map $\phi\mapsto p(.,\phi)$ is
of class $C^M(\T^n;S^m_N)$, one has $p_k\in S^m_N$ and 
$$
\norma{p_k}_{S^m_N}\leq\frac{\norma{p}_{C^M (\T^n;S^m_N)
}}{|k|^M}\ ,\quad k\not=0\ .
$$
Thus, defining 
\begin{equation}
\label{chi.22}
\chi(E,\phi,\omega):=\sum_{k\not=0}\frac{
  p_k(E)e^{\im k\cdot\phi} }{\im\omega\cdot k}\ ,
\end{equation}
for any $M_1<M-\tau-n$, one has 
$$
\norma{\chi}_{C^{M_1}(\T^n,S^m_N)}\leq \norma{p}_{C^M (\T^n;S^m_N)
}\sum_{k\not=0}\frac{|k|^{M_1+\tau}}{\gamma|k|^M}\ ,
$$
which is convergent. From the arbitrariness of $M$ it follows that
also $M_1$ can be chosen arbitrarily and therefore, for fixed
$\omega\in\Omega_\gamma$, the symbol $\chi\in \cS^m$. Furthermore,
since
\begin{equation}
\label{deromega}
\frac{\partial^n}{\partial \omega_j^n}\frac{1}{\im \omega\cdot
  k}=\frac{(-\im)^n n!k_j^n}{(-\im\omega\cdot k)^n}\ ,
\end{equation}
and similarly for the other derivatives,
one has that $\forall \rho$, the symbol $\chi\in
\widetilde{Lip}_\rho^m(\Omega_{0\gamma})$. 

Exploiting this remark it is easy to obtain the conclusion also for
the case of $p$ which depends on $\omega$ in a Whitney smooth
way.\qed 

In order to solve the equation \eqref{12.1} we define now the set 
\begin{equation}
\label{tigamma}
{\Omega_{1\gamma}}:=\left\{\omega\in\Omega\ :\ \left|\omega\cdot
k+k_0|\geq\frac{\gamma}{1+|k|^\tau}\right|\ ,\ (k_0,k)\in\Z^{n+1}-\left\{0\right\} \right\}\ .
\end{equation}

\begin{lemma}
\label{lemchi1}
Let $p\in Lip_\rho^{m}(\Omega_{1\gamma})$, then there exists a
solution $\chi\in {Lip}_\rho^{m}(\Omega_{1\gamma}) $ of
\eqref{12.p}.  Furthermore, in this case $\overline{\left\langle p\right\rangle }\in
\widetilde{Lip}_\rho^{m}(\Omega_{0\gamma}) $.
\end{lemma}

\proof Following \cite{Bam97}, we prove that the solution of the
homological equation \eqref{12.1} is given by
$$ \chi(x,\xi,\phi):=\sum_{k\in\Z^n}\chi_k(x,\xi)e^{\ir k\cdot \phi}\ ,
$$
where
\begin{align}
\label{ch.3}
\chi_0=\frac{1}{T}\int_0^{T}t(\overline p-\overline{\langle
  p\rangle})\circ \Phi^t_{h_0}dt
\\
\label{ch.4}
\chi_k(x,\xi)=\frac{1}{e^{\ir\omega\cdot k
    T}-1}\int_0^{T}e^{\ir \omega\cdot k
  t}p_k(\Phi^t_{h_0}(x,\xi))dt \ ,
\end{align}
and $p_k$ are the functions defined by \eqref{four}. To this end, 
consider the Eq. \eqref{12.1} and take first its $k$-th Fourier
coefficient (in $\phi$), we thus get 
\begin{equation}
\label{s.3}
\left\{\chi_k;h_0\right\}+\ir\omega\cdot k\chi_k=p_k-\delta_{k,0}p_0\ .
\end{equation}
For $k=0$ it reduces to \eqref{12} and thus we have already studied it.
For $k\not=0$ the equation \eqref{s.3} is the value at $t=0$ of the equation
\begin{equation}
\label{s.4}
\frac{d}{dt}\chi_k\circ\Phi^t_{h_0}+\ir\omega\cdot k
\chi_k\circ\Phi^t_{h_0}=p_k \circ\Phi^t_{h_0}\ ;
\end{equation}
denoting the r.h.s. by $p_k(t)$ we solve such an equation as an
ordinary differential equation for $\chi_k(t)$. The general solution
is given by
$$
a_ke^{-\ir\omega\cdot k t}+\int_0^te^{\ir\omega\cdot k s}p_k(s)ds\ .
$$
The value of the constant is determined by the requirement that the
solution must be periodic of period $T$. Thus one gets the formula
\eqref{ch.4}. Then it is immediate to use the Diophantine condition in
\eqref{tigamma} in order to estimate $\chi$ and its derivatives (both
in $\phi$ and in $\omega$). 
\qed

\subsection{The smoothing theorem}\label{smot}

We are now ready to state and prove the main result of the section.

\begin{theorem}
\label{smoothing}
Fix $\gamma>0$ small, $\rho>2$ and an arbitrary $\kappa>0$. Assume
\begin{align}
\label{sti.1}
\tilde \beta<l+1
\end{align} 
then there exists a (finite) sequence of symbols $\chi_1,...,\chi_N$
with $\chi_j\in Lip_\rho^{m_j}(\Omega_{0\gamma})$,
$m_{j} \geq m_{j+1}$ $\forall j$, s.t., defining 
\begin{equation}
\label{Xj}
X_j:=\chi^w_j(x,D_x,\omega t)\ ,\quad \omega\in\Omega_{0\gamma}\ ,
\end{equation}
such operators are selfadjoint and the transformation
\begin{equation}
\label{transfe}
\psi=e^{-\im \epsilon X_1(\omega t)}....e^{-\im\epsilon X_N(\omega t)}\varphi\ ,
\end{equation}
transforms $H_\epsilon(\omega t)$ (c.f. \eqref{H}) into a
pseudodifferential operator $H^{(reg)}=A_0+\epsilon R$, where 
\begin{equation}
\label{a0}
A_0:=H_0+\epsilon Z+\epsilon\tilde Z
\end{equation}
has symbol
\begin{equation}
\label{norfor}
h_0+\epsilon z(h_0)+\epsilon\tilde z(h_0,\omega)\ ,
\end{equation}
where $z\in \widetilde{S}^{\tilde \beta}$ is a function of $h_0$ independent of $\omega$; $\tilde z\in \widetilde{Lip}_\rho^{2\tilde\beta-l-1
}(\Omega_{0\gamma})$ is an $\omega$ dependent function of $h_0$
 and $R\in \cL ip_\rho^{-\kappa}(\Omega_{0\gamma})$ depends on
$(\phi,\omega)$. 
In the case $l=1$ the set $\Omega_{0\gamma}$ must be substituted by
the set $\Omega_{1\gamma}$.
\end{theorem}

\begin{remark}
\label{kappa}
Actually in order to develop the KAM part of the proof of Theorem
\ref{m.1} we only need the existence of a positive $\kappa$ s.t. the
above results hold.
\end{remark}

\noindent 
{\it Proof of Theorem \ref{smoothing} in the case $l>1$}. We only
study $h_0+\epsilon W_0$. We will transform it through a unitary (in
$L^2$) operator leaving invariant the spaces $\cH^s$; therefore under
such transformations
$W_{\infty}$ remains a smoothing operator of arbitrary order. 

Consider $h_0+\epsilon W_0$; we transform it using
the operator $X_1$ with symbol $\chi_1$ obtained by solving the
homological equation \eqref{12} with $p=W$, so that $\chi_1\in
\cS^{\beta-l+1}$, so that by Lemma \ref{MR} the corresponding Weyl
operator is selfadjoint provided
\begin{equation}
\label{onbeta}
\beta-l+1 \leq l+1\ \iff\ \beta\leq 2l
\end{equation}
and Lemma \ref{egorov} applies provided the inequality is strict. Then
the symbol of the transformed Hamiltonian is given by
\begin{align}
\label{h1.a}
h^{(1)}&:=h_0+\epsilon(\langle
W_0\rangle-W_0)+\epsilon^2\cS^{\beta-(l+1-m_1)}
+\epsilon \cS^{\beta-2l-2}
\\
&+\epsilon W_0+\epsilon^2  \left\{ W_0;\chi_1
\right\}+\epsilon^2 \cS^{\beta-2(l+1-m_1)}
\\
&-\epsilon \dot \chi_1+\epsilon^2 \cS^{2(\beta-l+1)-(l+1)}
\\
\label{h1.b}
&=h_0+\epsilon\langle W_0\rangle-\epsilon\dot\chi_1 +\epsilon P_1\ ,
\end{align}
with $P_1\in
\cS^{\beta-(2l-\beta)}$ and $m_1=\beta-l+1$.

Consider first the case where $\langle W\rangle\equiv 0$, which
implies $\langle W_0\rangle=0$. In this case we iterate the procedure
with $P_1':=-\dot \chi_1+ P_1\in \cS^{\beta_1}$ in place of $W_0$ and
$\beta_1:=max\left\{\beta-l+1;2\beta-2l\right\}$. Remark that
$$
\langle \dot\chi_1\rangle\equiv 0\ ,
$$
so that, after the second transformation generated by
$\chi_2\in S^{\beta_1-l+1}$,
one gets a Hamiltonian of the
form 
\begin{equation}
\label{tildeh1}
 h^{(1_2)}=h_0+\epsilon\langle P_1\rangle-\epsilon \dot
\chi_2+\epsilon\cS^{\beta_1-(2l-\beta_1)}\ .
\end{equation}
Then, if $\beta_1-l+1>2\beta-2l$ we iterate again until we get 
$$
\tilde h^{(1)}=h_0+\epsilon\langle P_1\rangle+\epsilon S^{\beta_2}\ ,
$$
with some $\beta_2< 2\beta-2l$.

In both cases we thus get (maybe after the second group of transformations) a
Hamiltonian of the form 
\begin{equation}
\label{htra1}
h^{(1')}:=h_0+\epsilon f(h_0,\omega t)+\epsilon P_2\ ,
\end{equation}
with $f(h_0,.)\in \cS^{\tilde \beta}$ and $P_2\in\cS^{\tilde\beta_1}$,
with 
\begin{align*}
\tilde \beta_1&<\tilde \beta
\end{align*}

We now continue by eliminating the time dependence from $f$. To this
end we take $\chi_3$ to be the solution of Eq. \eqref{12.p} with
$p=f(h_0)$, so that $\chi_3\in \widetilde{Lip}_\rho^{\tilde
  \beta}(\Omega_{0\gamma})$. Provided $$\tilde \beta<l+1\ ,$$ 
one gets that the corresponding Weyl
operator is selfadjoint and the quantum lie transform it generates,
transoms symbols into symbols. Then the symbol of the transformed
Hamiltonian takes the form
\begin{align*}
h^{(2)}=&h_0+\epsilon L ip_\rho^{2l+\tilde \beta-3(l+1)}
\\
&-\epsilon\dot \chi_3+\epsilon^2 L ip_\rho^{2\tilde \beta-3(l+1)}
\\
&+\epsilon f(h_0)+\epsilon^2 L ip_\rho^{2\tilde \beta-3(l+1)}
\\
&+\epsilon P_2+\epsilon^2 L ip_\rho^{\tilde\beta_1+\tilde\beta-(l+1)}
\\
=& h_0+\epsilon\overline{f(h_0)}+\epsilon P_2+\epsilon
L ip_\rho^{\tilde\beta_1+\tilde\beta-(l+1)}
\end{align*}
where all the functions are defined on $\Omega_{0\gamma}$. In particular
the perturbation $L ip_\rho^{\tilde\beta_1+\tilde\beta-(l+1)}$ is the
lowest order term with a nontrivial dependence on $\omega$.

Denote now
$$h_1:=h_0+\epsilon f(h_0)
$$ 
so that the Hamiltonian takes the form
$$
h_1+\epsilon P_2+\epsilon
L ip_\rho^{\tilde\beta_1+\tilde\beta-(l+1)}\ ;
$$ we can now iterate the above construction (with $h_1$ in place of
$h_0$, thus exploiting the homological equation \eqref{12.3}) until
we get a Hamiltonian of the form
\begin{equation}
\label{h2}
h^{(2')}:=h_2+\epsilon P_2'\ ,\quad P_2'\in L
ip_\rho^{\beta_2}\ ,\quad h_2:=h_1+\epsilon f_1(h_0)\ ,\quad
\beta_2:=\tilde\beta_1+\tilde\beta-(l+1)
\end{equation}
We are now in the position of concluding the proof of the theorem. We
proceed in a quite explicit way. First we construct $\chi_4$ by
solving \eqref{12.3} (with $h_2$ in place of $h_1$ and $P_2'$ in
place of $p$), we transform the Hamiltonian getting
\begin{equation}
\label{h3}
h_2+\epsilon\langle P_2' \rangle+  L
ip_\rho^{\beta_2-a_1}\ ,\quad a_1:=\min\left\{l-1;2l-\beta_2\right\}\ ,
\end{equation}
We construct now $\chi_5\in L ip_\rho^{\beta_2}$ by solving the
homological equation \eqref{12.p} and transforming the Hamiltonian
we get
\begin{equation}
\label{h3.1}
h^{(3)}=h_3+\epsilon L ip_\rho^{\beta_2-a_1}\ ,\quad
h_3:=h_2+\epsilon\overline{\langle P_2'\rangle} \ ;
\end{equation}
we remark that the correction in $h_3$ has a non trivial dependence on
$\omega$. 

We also remark that the gain in the order of the remainder does not decreases
as $\beta_2$ decreases. So one can iterate the construction lowering by
a finite quantity at each step the order of the perturbation. In this
way, by a finite number of step one gets the order $-\kappa$.
Finally remark that in the considered range of the parameters the
condition $\tilde\beta<l+1$ implies also condition \eqref{onbeta}.\qed

\vskip10pt

\noindent
{\it Proof of Theorem \ref{smoothing} in the case $l=1$.} First remark
that in this case the condition $\tilde \beta<l+1$ is equivalent to
$\beta<l+1$. We prove
that for any positive $N$ there exists
$\left\{\chi_j\right\}_{j=1}^N$, $\chi_j\in L
ip^{\beta-(j-1)(2-\beta)}(\Omega_{1\gamma})$
s.t. the symbol $h^{(N)}$ of the Hamiltonian obtained after the transformation
$\psi=e^{{\im}\epsilon X_N}...e^{{\im}\epsilon X_1}\phi$ has the structure
\begin{equation}
\label{h}
h^{(N)}= h_0+\epsilon\overline{\langle
W\rangle} +\epsilon^2 \widetilde{ z^{(N)}}+\epsilon^N r_N\ ,
\end{equation}
with $\widetilde{z^{(N)}}= \widetilde{z^{(N)}}(h_0,\omega)$,
$\widetilde{z^{(N)}}\in \widetilde S^{\beta-1}$ and
\begin{align}
\label{l=1.2}
r_N\in\cL ip^{\beta-N(2-\beta)}(\Omega_{1\gamma})\ .
\end{align}
We prove this by induction. Of course it is true for $N=0$ with
$r_0=W-\overline{\langle W\rangle}$. Assume it is true for $N$. We
transform now $h^{(N)}$ using 
$\epsilon^{N+1}\chi_{N+1}\in L
ip_\rho^{\beta_N}$, $\beta_N:=\beta-N(2-\beta)$,
which solves \eqref{12.1} with $p=r_N$. Remarking that in this case,
for any symbol $f$, one has 
$$
\left\{h_0,f\right\}^q=\left\{h_0,f\right\}\ ,
$$
it follows
\begin{align*}
h^{(N+1)}&=h_0+\epsilon^N\left\{h_0,\chi_{N+1}\right\}+\epsilon^{2N}\frac{1}{2}\left\{\left\{h_0;\chi_{N+1}\right\};\chi_{N+1}\right\}+l.o.t.
\\ &+\epsilon\overline{\langle W\rangle} +\epsilon^2 \widetilde{
  z^{(N)}}+\epsilon^{N+1} \left\{\overline{\langle
  W\rangle};\chi_{N+1}\right\}+l.o.t.
\\ &+\epsilon^Nr_N+\epsilon^{2N}\left\{r_N;\chi_{N+1}\right\}+l.o.t.
\\ &-\epsilon^N\dot\chi_{N+1}-\frac{1}{2}\epsilon^{2N}\left\{\dot\chi_{N+1};\chi_{N+1}\right\}
+l.o.t.  
\\ &=h_0+\epsilon\overline{\langle W\rangle} +\epsilon^2
\widetilde{ z^{(N+1)}}+\epsilon^{2N} L
ip_\rho^{2\beta_N-2}+\epsilon^{N+1}L ip_\rho^{\beta+\beta_N-2}\ ,
\end{align*}
where we put
$$
\widetilde{ z^{(N+1)}}:=\widetilde{ z^{(N)}}+\epsilon^{N-1}
\overline{\langle r_N\rangle}\ .
$$
\qed

\section{Diagonalization of the time independent part}\label{prepa}

In this section we diagonalize the operator 
\begin{equation}
\label{a0.1}
A_0:=H_0+\epsilon Z+\epsilon\tilde Z
\end{equation}
associated to the time independent part of the Hamiltonian:
\begin{equation*}
h_0+\epsilon z(h_0)+\epsilon\tilde z(h_0,\omega)\ ,
\end{equation*}
First write it in the basis $\be j$ of the normalized eigenvectors of
$H_0$ and fix a positive $s$ identifying the order of the space
$\cH^s$ in which we will control the norm of the operators and a
positive $\rho$ larger then 2 controlling the smoothness in $\omega$
of the various objects. 

We will denote $\Delta f:=f(\omega')-f(\omega)$ and
$\Delta\omega:=\omega'-\omega$.

\begin{lemma}
\label{diaga0}
There exists a positive $\epsilon_*$ s.t., if
$|\epsilon|< \epsilon_*$ then there exists a unitary (in
$L^2$) operator $U_1$, Whithey smooth in $\omega$, with 
\begin{equation}
\label{U1}
\norma{U_1-\uno}_{
  Lip_{\rho}(\Omega_{0\gamma};B(\cH^{s-\delta};\cH^{s})) }\sleq
\epsilon\ ,
\end{equation}
and $\delta:=\tilde \beta-(l+1)$, s.t.
\begin{equation}
\label{U2}
U_1^*A_0U_1=A^{(0)}\ ,
\end{equation}
where
\begin{align}
\label{A}
A^{(0)}&:=\diag(\lambda_j^{(0)})\ ,
\end{align}
with $\lambda_j^{(0)}$ given by
\begin{equation}
\label{diaga.6}
\lambda_j^{(0)}=\lambda_j^v+\epsilon z(\lambda_j^v)+\epsilon \tilde
z(\lambda_j^v,\omega) +\epsilon\frac{\nu_j(\omega)}{j^{\delta}}\ ;
\end{equation} 
and $\nu_j(\omega)$ Whitney smooth functions which fulfill
\begin{align}
\label{diaga.200}
\left|\nu_j(\omega)\right|\sleq 1\ ,
\\
\label{diaga.201}
\left|\frac{\Delta\nu_j(\omega)}{\Delta\omega}\right|\sleq 1
\end{align}
uniformly on $\Omega_{0\gamma}$ (or on  $\Omega_{1\gamma}$) and in $j$.
\end{lemma}
\proof Denote by $Z+\tilde Z$ the Weyl quantization of $z(h_0)+\tilde
z(h_0,\omega)$, then, from functional calculus one has 
that 
$$
R_a:=z(H_0)+\tilde z(H_0,\omega)-Z-\tilde Z\in Lip_\rho^{\tilde
  \beta-(l+1)}(\Omega_{0\gamma})\ ; 
$$
since $l+1>\tilde \beta$, the operator $R$ is smoothing of order
$\delta=l+1-\tilde \beta$. So we rewrite
\begin{equation}
\label{diaga.10}
A_0=\Lambda+\epsilon R_a\ ,\quad \Lambda:= H_0+\epsilon
z(H_0)+\epsilon\tilde z(H_0,\omega)\ .
\end{equation}
Then we diagonalize the system by a series of transformations which
are constructed in a way similar to the transformations that we will
use in sect. \ref{lemmaKAM} to prove Theorem \ref{analyticKAM} in
order to develop the KAM part of the proof. Here the situation is much
simpler since this procedure does not involve small denominators. In
order to develop the procedure we need to control the differences
between the eigenvalues. Denote
\begin{equation}
\label{diaga.400}
\lambda'_j:=\lambda_j^v+\epsilon z(\lambda_j^v)+\epsilon\tilde
z(\lambda^v_j,\omega)\ , 
\end{equation}
then we have to estimate from below $|\lambda'_i-\lambda'_j|$. To this
end consider first $z(\lambda_i^v)-z(\lambda_j^v)$, $i>j$. From the
mean value theorem there exists $\bar E\in(\lambda^v_j,\lambda^v_i)$
s.t.
\begin{equation}
\nonumber
|z(\lambda_i^v)-z(\lambda_j^v)|=\left|\frac{\partial z}{\partial E}(\bar
E)\right| |\lambda_i^v-\lambda_j^v|\sgeq
(\lambda_j^v)^{\frac{\tilde\beta}{l+1}-1}|\lambda_i^v-\lambda_j^v|\sgeq
|\lambda_i^v-\lambda_j^v|\ ,
\end{equation}
so that (repeating the argument for $\tilde z$), one has
\begin{equation}
\label{diaga.101}
|\lambda_i'-\lambda_j'|\sgeq |\lambda_i^v-\lambda_j^v|-\epsilon
|\lambda_i^v-\lambda_j^v|\sgeq |\lambda_i^v-\lambda_j^v|\sgeq
\left|i^d-j^d\right| \ .
\end{equation}
Define now an operator $X$ with matrix
elements  
$$
X_{ij}:=-\ir\frac{R_{a,ij}}{\lambda'_i-\lambda'_j}\ ,\quad
i\not= j
$$ so that, $-\ir[\Lambda;X]=\widebar R_a$, with $ \widebar
R_a=\diag(R_{a,ii})-R_a$. By Lemma \ref{pos92.c}, $X$ has the same
boundedness properties of $R_a$, so it is smoothing of order
$\delta$. Furthermore its norm is estimated by
$$
\norma{X}\sleq \norma{F}
$$ where the norm is the norm in
$Lip_\rho(\Omega_{0\gamma};B(\cH^{s-\delta};\cH^s))$, and the constant
depends only on the indexes of the norm and on the constant in the
inequality \eqref{diaga.101}. In this proof we will use only such
norm.

It follows
from lemma \ref{qliest} that the series defining $Lie_{\epsilon X} R_a$ is
convergent and one has
$$
\norma{Lie_{\epsilon X}F-F} \sleq\epsilon\norma{X}\norma{F}\sleq \epsilon \norma{F}^2\ .
$$
Furthermore exploiting the definition of $X$, one has 
\begin{align}
\label{diaga.50}
Lie_{\epsilon X}\Lambda =\Lambda +\epsilon \widebar R_a+\sum_{k\geq
  2}\epsilon^{k}\frac{\Lambda_k}{k!}\ ,
\\
\Lambda_1:=\widebar R_a\ ,\quad \Lambda_k=-\im [\Lambda_{k-1}; X ]\ ,
\end{align}
so that 
\begin{equation}
\label{trasf3}
Lie_{\epsilon X}(\Lambda+\epsilon R_a)=
\Lambda^{(1)}+\epsilon^2 R^{(1)}
\end{equation}
where $\Lambda^{(1)}:=\Lambda+\epsilon \diag (R_{a,ii})$ and $R^{(1)}$
is a suitable operator fulfilling
\begin{equation}
\label{diaga.102}
\norma{R^{(1)}}\sleq \norma{R_a}^2\ ,
\end{equation}
again with a constant which depends only on the indexes of the norm
and on the constant in the inequality \eqref{diaga.101}.

It is easy to see that the eigenvalues of $\Lambda^{(1)}$ fulfill
again inequality \eqref{diaga.101} with a constant which is decreased
by $O(\epsilon)$, so that one can iterate the argument and get the
existence of the operator $U_1$ claimed in the statement. 
\qed

We now study the properties of the eigenvalues \eqref{diaga.6}. Before
doing that, it is useful to introduce a few notations.
First we denote 
\begin{equation}
\label{diaga.210}
\inter m:=\max\left\{1;|m|\right\}\ .
\end{equation}
Then, given a closed set $\tilde\Omega\subset\Omega$, consider a
sequence $\lambda=\{\lambda_j(\omega)\}_{j\geq 1}$ of functions of $\omega$
defined on $\tilde\Omega$. We denote
\begin{align}
\label{sijk}
s_{ijk}(\lambda,\omega):=\lambda_i-\lambda_j+\omega\cdot k\ ,
\\
\label{rijk}
\cR_{ijk}(\lambda,\alpha):=\left\{\omega\in\tilde\Omega\ :\ \left|s_{ijk}(\lambda,\omega)\right|<\alpha\sidjd\right\}
\end{align}

The next lemma gives the properties of the eigenvalues.
We emphasize that in its proof we use the property that both $z$ and
$\tilde z$ are Whithney smooth in the frequencies.

Moreover, in the case $d=1$, we exploit the fact that $\delta>0$,
which is implied by $\tilde\beta<l+1=2$ (strictly). In the case $d>1$
this is not needed.

\begin{lemma}
\label{diaga.300}
There exists $\epsilon_*>0$ and $\tau>0$, s.t., for any
$|\epsilon|<\epsilon_*$ there exists $a>0$ and a closed set
$\Omega^{(0)}_\gamma\subset \Omega_{0\gamma}$ or
$\Omega^{(0)}_\gamma\subset \Omega_{1\gamma}$ with the following
properties
\begin{equation}
\label{diaga.301}
\left|\Omega-\Omega^{(0)}_\gamma\right|\sleq \gamma^a
\end{equation}
For any $\omega\in \Omega^{(0)}_\gamma$ the following inequalities
hold
\begin{align}
\label{diaga.1}
\left|\lambda_j^{(0)}-\lambda_j^v\right|\sleq j^{\frac{\tilde
    \beta}{l+1}}\ ,
\\
\label{diaga.2}
\left|\lambda_i^{(0)}-\lambda_j^{(0)}\right|\sgeq
\left|i^d-j^d\right|\ ,
\\
\label{diaga.3}
\left|\frac{\Delta(\lambda_i^{(0)}-\lambda_j^{(0)})}{\Delta\omega}
\right| 
\sleq
\epsilon|i^d-j^d| \ .
\\
\label{diaga.303}
\left|\lambda_i^{(0)}-\lambda_j^{(0)}+\omega\cdot k\right|\geq
\frac{\gamma\sidjd}{1+|k|^\tau}\ ,\quad
\left|i-j\right|+\left|k\right|\not=0\ .
\end{align}.
\end{lemma}
\begin{remark}
\label{a}
In the case $d>1$ one can choose $a=1$ and $\tau> n+2/(d-1)$. In the
case $d=1$ one can also compute such numbers, but they are more
complicated. 
\end{remark}

\proof Eqs. \eqref{diaga.1} and \eqref{diaga.2} immediately follow
from the previous proof. To get \eqref{diaga.3} compute
\begin{equation}
\label{delta.5}
\Delta(\lambda_i-\lambda_j)=\epsilon\Delta[\tilde
  z(\lambda_i^0)-\tilde z(\lambda_j^0)]+\epsilon\frac{\Delta
  \nu_i}{i^{\delta}}-\epsilon\frac{\Delta \nu_j}{j^{\delta}}\ ,
\end{equation}
To estimate the first term we use the mean value theorem (for Whitney
smooth function); to simplify the notation we
denote $$(\omega,\omega')=\{\nu\in \Omega_{0\gamma}\ :\ \exists
t\in(0,1)\ with\ \nu=t\omega+(1-t)\omega' \}\ .$$ So we have
\begin{align*}
\left|\Delta[\tilde z(\lambda_i^v)-\tilde z(\lambda_j^v)]\right|\leq
\sup_{\nu\in(\omega,\omega')}\left|(\omega'-\omega)\cdot
\frac{\partial}{\partial\omega}\left(\tilde z(\lambda_i^v,\nu)-\tilde
z(\lambda_j^v,\nu)\right)\right| 
\\ 
=
\sup_{\nu\in(\omega,\omega')}
\left|(\omega'-\omega)\cdot \left[\frac{\partial \tilde
    z}{\partial\omega}(\lambda_i^v,\nu)-\frac{\partial \tilde
    z}{\partial\omega}(\lambda_j^v,\nu)\right]\right|\leq
\left|\Delta\omega\right| \left|\lambda_i^v-\lambda_j^v\right|
\sup_{\nu\in(\omega,\omega')}
\sup_{\lambda\in(\lambda_j^v,\lambda_i^v)}\left|\frac{\partial^2\tilde
  z}{\partial\lambda\partial\omega}(\lambda,\nu)\right| 
\\ 
\sleq
\left|\Delta\omega\right| \left|\lambda_i^v-\lambda_j^v\right|\ .
\end{align*}
Adding the estimate of the other two terms one gets \eqref{diaga.3}. 

We come to \eqref{diaga.303}.  Define 
$$
\Omega^{(0)}_\gamma:=\Omega_{0\gamma}-\bigcup_{ijk}
\cR_{ijk}(\lambda^{(0)},\gamma/(1+|k|^{\tau}))\ .
$$

In order to estimate the above set we separate the case $d>1$ and the
case $d=1$. Consider first $d>1$; then by Lemma \ref{meas.4}, the
measure of $\cR_{ijk}(\lambda^{(0)},\gamma/(1+|k|^{\tau}))$ is estimated by
\eqref{meas.23} with $\alpha=\gamma/(1+|k|^{\tau})$. We fix $k$ and
estimate the cardinality of the $i,j$'s such that the set $\cR_{ijk}$ is
not empty. By \eqref{notempty}, exploiting the fact that
$$
\idjd\sgeq
\left(i^{d-1}+j^{d-1}\right)\left|i-j\right|\geq(i^{d-1}+j^{d-1})\ , 
$$
such a cardinality is estimated by $|k|^{2/(d-1)}$, so we have
$$
\left|\bigcup_{ijk}
\cR_{ijk}(\lambda^{(0)},\gamma/(1+|k|^{\tau}))\right|\sleq
\sum_{k\in\Z^n}\frac{\gamma |k|^{2/(d-1)}}{1+|k|^\tau}\sleq \gamma\ ,
$$
which concludes the proof in the case $d>1$. 

The case $d=1$ is slightly more complicated. In this case we have
$\lambda_j^v= j+\frac{1}{2}$, so that
the cardinalty above is infinite.

First we write $i=j+m$, so that one has
\begin{equation}
\lambda_i^{(0)}-\lambda^{(0)}_j=m+\epsilon
(f(j+m+1/2)-f(j+1/2))+\epsilon
\left(\frac{\nu_{j+m}}{(j+m)^\delta}-\frac{\nu_{j}}{j^\delta}\right)
\ ,
\end{equation}
where $f=z+\tilde z$. Now, by the mean value theorem, there exists
$\bar E\in (j+1/2,j+m+1/2)$ s.t.
$$
|f(j+m+1/2)-f(j+1/2)|= |f'(\bar E)m|\leq \frac{2}{j^{1-\tilde
    \beta/2}} m\ .
$$
Let $C$ be the constant in \eqref{diaga.200}, and define 
$$
\delta_{j}:=\frac{2}{j^{1-\tilde
    \beta/2}}+2\frac{C}{j^{\delta}}
$$
so that $|s_{ijk}-m-\omega\cdot k|\leq m\delta_{j}$. 
Define now the sets 
\begin{equation}
\label{qijk}
\cQ_{mjk}=\left\{\omega\in\Omega_{1\gamma}\ :\ \left|m+\omega\cdot k
\right|<\frac{\gamma m}{1+|k|^\tau}+m\delta_j   \right\}\ ,
\end{equation}
and remark that $\cR_{ijk}\subset \cQ_{mjk}$ and also
$\cQ_{mjk}\subset \cQ_{mj'k}$ if $j>j'$. Exploiting this remark we
take some $j_*$, fix $k$ and proceed as follows
\begin{equation}
\label{meas.10}
\bigcup_{ij}\cR_{ijk}(\lambda^{(0)},\alpha)\subset
\left(\bigcup_{i-j=m, j<j_*}\cR_{ijk}\right)\cup
\left(\bigcup_{m}\cQ_{mj_*k} \right)\ ;
\end{equation}

then (by Lemma \ref{meas.4}), remarking that \eqref{notempty} implies
$|m|\sleq |k|$
 one has
$$
\left|\left(\bigcup_{i-j=m, j<j_*}\cR_{ijk}\right)\right|\sleq \alpha
j_*|k| \ .
$$ 
Furthermore, one has that the set $\cQ_{mjk}$ is just the set
$\cR_{ijk}$ with $\lambda_j=j$ and
$\alpha=\gamma/(1+|k|^\tau)+\delta_j$. It follows that
$$
\left|\bigcup_{m}\cQ_{mj_*k} \right|\sleq |k| \left(\frac{\gamma
  }{1+|k|^\tau}+\delta_{j_*} \right)\ ,
$$
and therefore the measure of \eqref{meas.10} is estimated by
$$
\left| k\right|\left(\frac{\gamma
  (j_*+1)}{1+|k|^{\tau}}+\delta_{j_*}\right)\simeq |k|\left(\frac{\gamma
  j_*}{1+|k|^{\tau}}+\frac{1}{j_*^{\tilde\delta}} \right)
$$
where $\tilde \delta=\min\left\{\delta;1-\tilde\beta/2\right\}$;
choosing
$j_*=\left(\frac{1+|k|^{\tau}}{\gamma}\right)^{1/(\tilde\delta+1)}$,
inserting in the above estimate and summing over $k$ one gets the thesis.
\qed

\begin{corollary}
\label{Hop}
The transformation $U_1$ transform $A_0+\epsilon R$ into 
\begin{equation}
\label{H0R}
A^{(0)}+\epsilon R_0\ ,
\end{equation}
where 
\begin{equation}
\label{HoR.1}
R_0:=U_1^{-1}RU_1\in \Lip_\rho(\Omega^{(0)}_\gamma;
C^\ell(\T^n;B(\cH^{s-\kappa};\cH^s)))\ , \quad \forall \ell\ .
 \end{equation}
\end{corollary}

\section{Analytic KAM theory}\label{lemmaKAM}

In this section we prove KAM theorem for analytic perturbations of
$A^{(0)}$. The procedure is essentially identical to the one developed
in \cite{BG01} (which is actually a small modification of
\cite{pos96}), except that we take here advantage of the fact that the
perturbation is smoothing, so everything is slightly simpler.

In the previous section we fixed a positive arbitrary $s$; now we also
fix a positive (large) $\kappa$, then we define the following norms of
operators and of operator valued functions of $\omega\in\tilde \Omega$
(with $\tilde \Omega$ closed), and of $\phi\in\T^n_r$. Here and below
$\T^n_{r}$ is the set of the angles belonging to the complexified
torus and fulfilling $|\Im \phi_j|<r$.

Let $F:\T^n_r\mapsto B(\cH^{s-\kappa},\cH^s)$, be an analytic map. We
define
\begin{equation}
\label{normas}
\norma{F}_r:=\sup_{\phi\in\T^n_r}\norma{F(\phi)}_{
  B(\cH^{s-\kappa},\cH^s)) }\ .
\end{equation}  
If $F$ depends also in a Lipschitz way on $\omega\in\widetilde
\Omega$, we still denote 
\begin{equation}
\label{normas.1}
\norma{F}_r:=\sup_{\omega\in\tilde
  \Omega}\sup_{\phi\in\T^n_r}\norma{F(\phi,\omega)}_{
  B(\cH^{s-\kappa},\cH^s)) }\ ,
\end{equation}  
and we define
\begin{equation}
\label{normalips}
\normaL{F}_r:=\sup_{\omega\not=\omega'\in\tilde
  \Omega}\frac{\norma{F(\omega)-F(\omega')}_r}{|\omega-\omega'|}\ .
\end{equation}

\begin{definition}
\label{def.anli}
An analytic map $F$ which is Lipschitz dependent on
$\omega\in\widetilde \Omega$ will be said to be Lipschitz analytic. 
\end{definition}

\subsection{Squaring the order of the perturbation}\label{squaring}

Consider a Hamiltonian of the form 
\begin{equation}
\label{H.a}
H=A+P(\omega t,\omega)\ ,\quad A=\diag(\lambda_j(\omega))\ ,\quad
\norma{P}_r,\normaL{P}_r\ll 1\ .
\end{equation}
We look for a selfadjoint operator $X=X(\omega t)$ with the property
that the transformation $T_X$ that it generates according to
Definition \ref{t.tras}
transforms $H$ into
\begin{equation}
\label{H+}
H^+=A^++P^+(\omega t,\omega)\ ,\quad A^+=\diag(\lambda^+_j(\omega))\ ,
\end{equation}
with $P^+$ having a size which is essentially the square of that of
$P$.

By Lemma \ref{T.1} one has 
\begin{align}
\label{squar.1}
H^+&= A
\\
\label{squar.2}
&-\ir[A;X]-\dot X+P
\\
\label{squar.3}
&+Lie_XA -\left(A-\ir[A;X]\right) 
\\
\label{squar.4}
&+Lie_XP-P 
\\
\label{squar.5}
&+Y_X+\dot X \ .
\end{align}
So, we look for an $X$  solving the ``quantum homological equation'':
\begin{equation}
\label{qhom}
-\ir[A;X]-\dot X+P=[P]\ ,\quad [P]:=\diag(P_{jj})\ ,
\end{equation}
in order to get the wanted result with
\begin{equation}
\label{qhom.1}
A^+=A+[P]\ ,\quad P^+=\eqref{squar.3}+\eqref{squar.4}+\eqref{squar.5}\ .
\end{equation}

\begin{lemma}
\label{solcohomoq}
Fix positive constants $\Gamma,\cK_1$. Assume that there exists a set
$\Omega_\gamma\subset \widetilde \Omega$ s.t. $\forall
\omega\in\Omega_\gamma$ one has
\begin{align}
\label{diofq}
\left|\lambda_i(\omega)-\lambda_j(\omega)+\omega\cdot
k\right|\geq\frac{\gamma\sidjd}{\ak}\ ,\quad \ijkno
\\
\label{diofq.2}
\left|\frac{\Delta(\lambda_i-\lambda_j)}{\Delta\omega}\right|\leq
\cK_{1\lambda} \idjd
\end{align}
with some 
\begin{equation}
\label{ite.20}
\gamma\geq \Gamma\ ,\quad \cK_{1\lambda}<\cK_1\ ,
\end{equation}
then Eq. \eqref{qhom} has an analytic Lipschitz solution $X$ defined on
$\Omega_\gamma$ and fulfilling 
\begin{align}
\label{estisolch.1}
\norma
X_{r-\sigma}\sleq\frac{\norma{P}_r}{\sigma^{n+\tau}}\ ,
\\
\label{estisolch}
\normaL
X_{r-\sigma}\sleq\frac{\normaL{P}_r}{\sigma^{n+\tau}}+
\frac{\norma{P}_r}{\sigma^{n+2\tau+1}}\ .
\end{align}
where the constants depend on $\Gamma$, $\tau$, $n$.
\end{lemma}
\proof The proof is standard. We insert it only for the sake of
completeness. Expanding the equation \eqref{qhom} in Fourier
series 
$$
X_k:=\frac{1}{(2\pi)^n}\int_{\T^n} X(\phi)e^{-\ir k\cdot \phi}d\phi
$$
one gets
$$
\ir[A,X_k]-\ir\omega\cdot k X_k= P_k-[P_k]\ ;
$$ 
taking the $ij$ element of matrix one gets that $X$ can be defined by
\begin{equation}
\label{s.c.1}
X_{kij}:=\frac{P_{kij}}{\ir(\lambda_i-\lambda_j+\omega\cdot k)}\ ;
\end{equation}
remark now that one has
\begin{equation}
\label{stiex}
\norma{P_k}_{B(\cH^{s-\kappa};\cH^{s})}\leq
\norma{P}_re^{-|k|r}\ ,\quad \frac{\norma{\Delta
    P_k}_{B(\cH^{s-\kappa};\cH^{s})}}{|\Delta\omega|}\leq
 \normaL{P}_re^{-|k|r}\ ,
\end{equation}
so that
$$
\left| X_{kij}\right|\sleq
\frac{|P_{kij}|(\ak)}{\gamma\sidjd}\ .
$$
Applying Lemma \ref{pos96}, one gets
\begin{equation}
\label{stipk}
\norma{X_k}_{r-\sigma}\sleq \frac{\norma{P_k}(\ak)}{\gamma}\sleq
\frac{\norma{P}_r(\ak)e^{-r|k|}}{\gamma}\ , 
\end{equation} 
and 
\begin{equation}
\label{norX}
\norma{X}_{r-\sigma} \sleq
\frac{\norma{P}_r}{\gamma}\frac{1}{\sigma^{n+\tau}}
\left[\sigma^n\sum_{k\in\Z^n}(1+\left|\sigma
  k\right|^\tau)e^{-\sigma|k|} \right]\ ,
\end{equation}
and remarking that the term in square bracket is just the Riemann sum
for the integral of the function $(1+|y|^{\tau})e^{-|y|}$ one gets
\eqref{estisolch.1}. 

To get \eqref{estisolch} write
\begin{equation}
\label{esich.3}
|\Delta X_{kij}|\leq\left|\Delta\frac{ P_{kij}}{s_{ijk}}\right|\leq
\left|\frac{\Delta
  P_{kij}}{s_{ijk}}\right|+\left|P_{kij}\right|\left|\Delta\frac{1}{s_{ijk}}\right|
\ ;
\end{equation}
The first addendum is estimated exactly as before. Concerning the
second one, one has
\begin{align}
\left|\Delta\frac{1}{s_{ijk}}\right|=\left|
\frac{s_{ijk}(\omega')-s_{ijk}(\omega)}{s_{ijk}(\omega')
  s_{ijk}(\omega)}\right|\leq
\frac{|k||\Delta\omega|+|\Delta(\lambda'_{i}-\lambda'_j) |}
     {s_{ijk}(\omega') s_{ijk}(\omega)} 
\\ 
\leq
     |\Delta\omega|(\ak)^2\frac{|k|+\cK_{1\lambda}\idjd}
     {\gamma^2\sidjd^2}
\\
\leq 
     |\Delta\omega|(\ak)^2(1+|k|)\idjd\frac{(1+\cK_{1\lambda})}
     {\gamma^2\sidjd^2}
\\
\leq |\Delta\omega|\frac{(\ak)^2(1+|k|)}
     {\gamma^2\sidjd}\ .
\end{align}
Then, proceeding exactly as in the previous case one gets the thesis.\qed

Exploiting the Lemmas \ref{reLie}, \ref{estia} and \ref{lemmay}, it is
immediate to get the following result: 
\begin{lemma}
\label{square}
Under the assumptions of Lemma \ref{solcohomoq}, there exists a
constant $C_*$ s.t., if
$$
\norma{X}_{r-\sigma}\leq C_*
$$
then one has 
\begin{equation}
\label{sqa.1}
\norma{P^+}_{r-\sigma}\sleq
\frac{\norma{P}_r}{\sigma^{2n+2\tau+1}}\norma P_r \ ,\quad
\normaL{P^+}_{r-\sigma}\sleq
\frac{\norma{P}_r}{\sigma^{2n+2\tau+1}}\left(\normaL
P_{r}+\frac{\norma P_r}{\sigma^{\tau+1}}\right) \ ,
\end{equation}
furthermore one has
\begin{equation}
\label{defor}
\norma{\uno-e^{-iX}}_{r-\sigma}\sleq 
\frac{\norma P_r}{\sigma^{n+\tau}}\ ,\quad \normaL{\uno-e^{-iX}}_{r-\sigma}\sleq 
\frac{\normaL{P}_r}{\sigma^{n+\tau}}+
\frac{\norma{P}_r}{\sigma^{n+2\tau+1}}\ .
\end{equation}
\end{lemma}

In order to be able to iterate the construction we still have to show
that the new eigenvalues also fulfill a Diophantine inequality.

\begin{lemma}
\label{autov}
Fix two constants $\cK_0$ and $\cK_1$ fulfilling 
\begin{equation}
\label{autov.300}
\cK_1\leq\frac{\cK_0}{8}\ .
\end{equation}
Assume $\kappa>1$ and $\tau>\left(n-\kappa^{-1}\right)/(1-\kappa^{-1})$.
Assume also that the eigenvalues $\lambda_j$ fulfill the following
estimates
\begin{align}
\label{autov.1}
\left|\lambda_i-\lambda_j\right|\geq \cK_{0\lambda}\idjd\ ,
\\
\label{autov.2}
\left|\frac{\Delta(\lambda_i-\lambda_j)}{\Delta\omega}\right|\leq
\cK_{1\lambda} \idjd\ ,
\end{align}
and \eqref{diofq}. Fix some $K$ fulfilling $1+K^\tau\leq \gamma/\norma
P_r$, then the eigenvalues $\lambda_j^+$ fulfill \eqref{autov.1} and
\eqref{autov.2} with new constants given by
\begin{align}
\label{autov.4}
\cK^+_{0\lambda}:=\cK_{0\lambda}-2 \norma{P}_r\ ,
\\
\label{autov.5}
\cK^+_{1\lambda}:=\cK_{1\lambda}+2 \normaL{P}_r\ ,
\end{align}
Assume furthermore that 
$$
\cK_{0\lambda}^+\geq \cK_0\ ,\quad \cK_{1\lambda}^+\leq \cK_1
$$
then there exists a measurable set $\Omega^+_{\gamma}$, and a positive
constant $d_1$, such that, for
any $\omega\in \Omega^+_{\gamma^+}$ one has 
\begin{equation}
\label{diofq.1}
\left|\lambda^+_i(\omega)-\lambda^+_j(\omega)+\omega\cdot
k\right|\geq\frac{\gamma^+\sidjd}{\ak}\ ,\quad\ijkno
\end{equation}
with 
\begin{equation}
\label{autov.301}
|\Omega_\gamma-\Omega_{\gamma^+}^+|\sleq \gamma^+ K^{-d_1}
\end{equation}
and 
\begin{equation}
\label{autov.6}
\gamma^+=\gamma-2(1+K^\tau)\norma P_r\ .
\end{equation}
The constant in \eqref{autov.301} depends on $\tau,$ $\kappa$,
$\cK_0,$ $\cK_1$ and $d_1$ 
\end{lemma}
\proof Denote $\epsilon_1:=\norma P_r$ and $\epsilon_2:=\normaL
P_r$. We have 
$$
\lambda_j^+=\lambda_j+\frac{\nu_j(\omega)}{j^\kappa}\ ,\quad
|\nu_j|\leq \epsilon_1\ ,\quad
\left|\frac{\Delta\nu_j}{\Delta\omega}\right| \leq \epsilon_2\ .
$$
Therefore \eqref{autov.4} and \eqref{autov.5} immediately
follow. Furthermore one has
\begin{align}
\label{autoval.100}
\left|\lambda^+_i(\omega)-\lambda^+_j(\omega)+\omega\cdot
k\right|\geq\left|\lambda_i(\omega)-\lambda_j(\omega)+\omega\cdot
k\right|-\frac{2\epsilon_1}{j^{\kappa}}\geq \frac{\sidjd}{\ak}\left(
\gamma-\frac{2\epsilon_1}{j^\kappa}\frac{\ak}{\sidjd}\right)\ ,
\end{align}
which is automatically larger then the r.h.s. of \eqref{diofq.1} if
\begin{equation}
\label{autoval.101}
\frac{\ak}{j^\kappa\sidjd}\leq 1+K^\tau\ .
\end{equation}
In turn this is automatic if $|k|\leq K$. So we consider now the case
$|k|>K$. In that case \eqref{autoval.101} is again automatic if
$j^\kappa>|k|^{\tau-1}\frac{\cK_{0\lambda}}{4}$ (exploiting
\eqref{notempty}). 

So, fix a value of $k$ with $|k|>K$ and consider the case
\begin{equation}
\label{autoval.102}
j\leq
|k|^{\frac{\tau-1}{\kappa}}\left(\frac{\cK_{0\lambda}}{4}\right)^\frac{1}{\kappa}
\ .
\end{equation}
Fix a value of $j$ fulfilling \eqref{autoval.102}, then (from
\eqref{notempty}) the set
$\cR_{ijk}^+:=\cR_{ijk}(\lambda^+,\gamma^+/(\ak))$ is not empty only
for a set of $i$'s which has at most a cardinality proportional to
$|k|$. We have (by Lemma \ref{meas.4})
\begin{equation}
\label{autoval.104}
\left|\bigcup_{i,j}\cR_{ijk}^+\right|\sleq
\sum_{i,j}\frac{1}{\cK_{0\lambda}^+} \frac{\gamma^+}{\ak}\sleq
\frac{1}{\cK_{0\lambda}^+} \frac{\gamma^+}{\ak} |k| |k|^{
  \frac{\tau-1}{\kappa}}\ ,
\end{equation}
where the sum is restricted to the $j$'s fulfilling
\eqref{autoval.102} and the $i$'s for which $\cR_{ijk}^+$ is not empty.
Summing over $k$ with $|k|>K$ one gets the result. \qed

\subsection{Iterative lemma and Analytic KAM theorem}\label{kamana}

We are now in the position of stating the iterative Lemma which is a
direct consequence of the results of the above subsection. Such a Lemma
yields the analytic KAM result that we need.

To start with take a positive $r$ and consider a quantum Hamiltonian
of the form
\begin{equation}
\label{H0}
H^{(0)}=A^{(0)}+P^{(0)}\ ,
\end{equation}
with $A^{(0)}$ given by \eqref{A} and $P^{(0)}$ an analytic Lipschitz
map fulfilling 
\begin{equation}
\label{ite.9}
\norma{P^{(0)}}_{r}\leq \epsilon_1^{(0)}\ ,\quad
\normaL{P^{(0)}}_{r}\leq \epsilon_2^{(0)}\ ,
\end{equation}
with some positive (small) $\epsilon_1^{(0)}$ and $\epsilon_2^{(0)}$.

The next lemma is a direct consequence of Lemmas \ref{square} and
\ref{autov} applied iteratively by taking $K=K^{(l)}$ as defined by
the first of \eqref{ite.8}. 

\begin{lemma}
\label{iterative}
Fix $0<\vartheta<1$, $\cK_0$, $\cK_1$ and $\Gamma$ with $\cK_0\geq
8\cK_1$ and define 
\begin{equation}
\label{ite.1}
\sigma_l:=\frac{(1-\vartheta) r}{2^l}\ ,\quad
r_l:=r-\sum_{i=1}^{l}\sigma_i\ .
\end{equation}
Then, there exist positive constants $d_2$, $d_3$ s.t., if one defines
iteratively (for $l\geq 0$) the sequences of constants
\begin{align}
\label{ite.1.1}
\epsilon_1^{(l+1)}\simeq\frac{(\epsilon_1^{(l)})^2}
        {(\sigma^{(l)})^{2n+2\tau+1}}  \ ,\quad
        \epsilon_2^{(l+1)}\simeq \frac{\epsilon_1^{(l)}}
        {(\sigma^{(l)})^{2n+2\tau+1}}
        \left(\epsilon_2^{(l)}+\frac{\epsilon_1^{(l)}}{
          (\sigma^{(l)})^{\tau+1}}\right) \ ,
\\
\label{ite.8}
K^{(l)}=(\epsilon_1^{(l)})^{-1/2\tau}\ ,\quad
\gamma^{(l+1)}=\gamma^{(l)}-(\epsilon_1^{(l)})^{d_2}\ ,\quad \delta^{(l)}\simeq \gamma^{(l)}   (\epsilon_1^{(l)})^{d_3}\ ,
 \\
\label{ite.91}
\cK^{(l+1)}_{0\lambda}=\cK^{(l)}_{0\lambda}-2\epsilon_1^{(l)}\ ,\quad
\cK^{(l+1)}_{1\lambda} =\cK^{(l)}_{1\lambda}+2\epsilon_1^{(l)}\ .
\end{align}
and for any $l\geq 0$ the following inequalities hold
\begin{align}
\label{ite.101}
\frac{\epsilon_1^{(l)}}
        {(\sigma^{(l)})^{n+\tau}}\sleq 1\ ,\quad \frac{\epsilon_2^{(l)}}
        {(\sigma^{(l-1)})^{n+\tau}}
        +\frac{\epsilon_1^{(l)}}{
          (\sigma^{(l)})^{n+2\tau+1}}\sleq 1 
\\
\label{ite.102}
\gamma^{(l)}\geq \Gamma\ ,\quad \cK_{0\lambda}^{(l)}\geq
\cK_0\ ,\quad \cK_{1\lambda}^{(l)}\leq \cK_1\ ,
\end{align}
then the following holds true: for any $l$ there exists a measurable
set $\Omega^{(l)}_{\gamma^{(l)}}$ and a Lipschitz analytic map $X^{(l)}$ defined on
$\Omega^{(l)}_{\gamma^{(l)}}$ with the property that $T_{X^{(l)}}$ is
well defined and one has
\begin{equation}
\label{ite.3}
T_{X^{(l)}}H^{(l-1)}=H^{(l)}=A^{(l)}+P^{(l)}\ ,\quad l\geq 1\ ,
\end{equation} 
with $A^{(l)}=\diag(\lambda_j^{(l)})$. Furthermore the following estimates
hold
\begin{align}
\label{ite.4}
\left|\Omega^{(l-1)}_{\gamma^{(l-1)}}-\Omega^{(l)}_{\gamma^{(l)}}\right|\leq
\delta^{(l)}\ , \\
\label{ite.5}
\norma{P^{(l)}}_{r_l}\leq \epsilon_1^{(l)}\ ,\quad
\normaL{P^{(l)}}_{r_l}\leq \epsilon_2^{(l)}\ ,
\\
\label{ite.6}
\left|\lambda_i^{(l)}-\lambda_j^{(l)}\right|\geq \cK^{(l)}_{0\lambda}
\idjd \ ,\quad \left|\frac{\Delta(\lambda_i^{(l)}-\lambda_j^{(l)}
  )}{\Delta\omega} \right|\leq \cK^{(l)}_{1\lambda}
\idjd \ ,
\\
\label{ite.7}
\left|\lambda_i^{(l)}-\lambda_j^{(l)}+\omega\cdot k
\right|\geq\frac{\gamma^{(l)} }{\ak}\sidjd\ ,\quad \ijkno
\\
\label{ite.11}
\norma{X^{(l)}}_{r_l}\sleq
\frac{\epsilon_1^{(l-1)}}
        {(\sigma^{(l-1)})^{n+\tau}}  \ ,\quad
\normaL{X^{(l)}}_{r_l}\sleq\frac{\epsilon_2^{(l-1)}}
        {(\sigma^{(l-1)})^{n+\tau}}
        +\frac{\epsilon_1^{(l-1)}}{
          (\sigma^{(l-1)})^{n+2\tau+1}} \ .
\\
\label{ite.103}
\norma{\uno-e^{-\ir X^{(l)}}}_{r_l}\sleq \norma {X^{(l)}}_{r_l} \ ,\quad
\normaL{\uno-e^{-\ir X^{(l)}}}_{r_l}\sleq \normaL{X^{(l)}}_{r_l}
\end{align}
\end{lemma}

\begin{theorem}
\label{analyticKAM}
Consider the quantum Hamiltonian \eqref{H0}, defined and Lipschitz on
a set $\Omega^{(0)}_{\gamma^{(0)}}$ s.t.
\begin{equation}
\label{an.1}
\left|\Omega-\Omega^{(0)}_{\gamma^{(0)}}\right|\leq \Upsilon^{(0)}
\end{equation}
with some positive $\Upsilon^{(0)}$. Fix positive numbers $\cK_0,$
$\cK_1,$ $\Gamma,$ $\tau,$ $\vartheta$ fulfilling
\begin{equation}
\label{an.400}
\tau>\frac{n-\kappa^{-1}}{1-\kappa^{-1}}\ ,\quad 0<\vartheta<1\ ,\quad
\cK_0>8\cK_1\ .
\end{equation}
Assume that, for some $0<r\leq 1$ and some positive $\varsigma$ one has
\begin{equation}
\label{hip}
\norma{P^{(0)}}_r\leq\varsigma r^b\ ,\quad
r^{\tau+1}\normaL{P^{(0)}}_r\leq\varsigma 
r^b\ ,\quad b:=2n+2\tau+1\ .
\end{equation}  
Then there exist positive constants
$\varsigma_*,C_\Gamma,C_0,C_\Omega, C_U$ s.t.,
if $|\varsigma|<\varsigma_*$ and the eigenvalues $\lambda_j^{(0)}$
fulfill 
\begin{align}
\label{an.2}
\left|\lambda^{(0)}_i-\lambda^{(0)}_j\right|\geq \cK^{(0)}_{0\lambda}
\idjd
\\
\label{an.3}
\left|\frac{\Delta(\lambda^{(0)}_i-\lambda^{(0)}_j)}{\Delta\omega}\right|
\leq \cK^{(0)}_{1\lambda}\idjd\ ,
\\
\label{an.4}
\left| \lambda^{(0)}_i-\lambda^{(0)}_j+\omega\cdot k \right|\geq
\frac{\gamma^{(0)}}{\ak}\sidjd\ ,\ \ijkno
\end{align}
with constants s.t 
\begin{align}
\label{an.5}
\cK^{(0)}_{0,\lambda}-C_0r^b\vsi >\cK_0\ ,\quad
\cK^{(0)}_{1\lambda}+C_0 r^{b+\tau+1}\vsi <\cK_1\ , \quad
\gamma^{(0)}-C_{\Gamma}(r^b\vsi)^{d_2}\geq \Gamma\ ,
\end{align}
 then there exists a measurable set $\Omega^{(\infty)}_{\gamma^{(\infty)}}$ and 
a Lipschitz analytic map $U$ defined on $\Omega^{(\infty)}_{\gamma^{(\infty)}}$, with
$U(\phi,\omega)$,  $L^2$ unitary, s.t.
the transformation
  $U(\omega 
  t,\omega)\psi'=\psi $ transforms the system \eqref{H0} into
\begin{equation}
\label{indep}
\ir\dot \psi'=A^{(\infty)}\psi'\ ,\quad
A^{(\infty)}:=\diag(\lambda^{(\infty)}_j)\ .
\end{equation}
Furthermore the following estimates hold
\begin{align}
\label{an.21}
\left|\lambda_i^{(\infty)}-\lambda_j^{(\infty)}\right|\geq (\cK^{(0)}_{0\lambda}-C_0r^b\vsi)
\idjd
\\
\label{an.31}
\left|\frac{\Delta(\lambda^{(\infty)}_i-\lambda^{(\infty)}_j)}{\Delta\omega}\right|
\leq (\cK^{(0)}_{1\lambda}+C_0  r^{b+\tau+1}\vsi )\idjd\ ,
\\
\label{an.41}
\left| \lambda^{(\infty)}_i-\lambda^{(\infty)}_j+\omega\cdot k \right|\geq
\frac{\gamma^{(0)}-C_{\Gamma}(r^b\vsi)^{d_2} }{\ak}\sidjd\ ,\ \ijkno
\\
\label{an.51}
|\Omega^{(0)}_{\gamma^{(0)}}-\Omega^{(\infty)}_{\gamma^{(\infty)}}|\leq \Upsilon^{(0)}(1+ C_{\Omega}
\gamma^{(0)}(\vsi r^b)^{d_3})\ ,
\\
\label{stitrakam}
\norma{\uno-U}_{\theta r}\leq C_U\vsi r^{b-(n+\tau)}\ ,\quad
r^{\tau+1}\normaL{\uno-U}_{\theta r}\leq C_U\vsi r^{b-(n+\tau)}\ .
\end{align}
\end{theorem}

\proof We apply Lemma \ref{iterative}. To this end we define
$\epsilon_1^{(l)}:=\vsi^{(l)} r^b$ and $\epsilon_2^{(l)}:=\vsi^{(l)}
r^{b-\tau-1}$ with $\vsi^{(0)}:=\vsi$. We fix $\vartheta=\frac{1}{2}$,
then all the constants \eqref{ite.1.1}-\eqref{ite.91} are defined by
the recursion. We first analyze \eqref{ite.1.1} which take the form
\begin{equation}
\label{stimerec}
\epsilon_1^{(l+1)}\sleq \frac{\epsilon_1^{(l)}2^{bl}}{
  r^b}\epsilon_1^{(l)} \ ,\quad \epsilon_2^{(l+1)}\sleq
\frac{\epsilon_1^{(l)}2^{bl}}{
  r^b}
\left(\epsilon_2^{(l)}+\frac{\epsilon_1^{(l)}2^{l(\tau+1)}}{  r^{\tau+1}} \right)\ ,
\end{equation}
which in turn can be reformulated in terms of $\vars l$:
\begin{align*}
\vars{l+1} r^b\sleq \frac{\vars lr^b2^{bl}}{r^b}\vars
lr^b\ \iff \vars{l+1}\sleq {2^{bl}}\left(\vars l\right)^2
\\ \vars{l+1} r^{b-\tau-1}\sleq \frac{\vars lr^b2^{bl}}{ r^b}
\vars l \left(\vars l r^{b-(\tau+1)}+\frac{\vars
  l2^{l(\tau+1)}r^b}{ r^{\tau+1}}\right)
\ \iff \vars{l+1}\sleq {2^{(b+\tau+1)l}}\left(\vars
l\right)^2\ ,
\end{align*}
which is solved, thanks to Lemma \ref{rico}, by defining
\begin{equation}
\label{ber}
\vars l:=\frac{1}{c_3
  2^{(b+\tau+1)l}}\left(2^{b+\tau+1}c_3\vsi
\right)^{2^{l}}\ ;\quad l\geq 1\ ,
\end{equation}
with $c_3$ the non written constant in the definition of the iterative
estimates. Then $\vars l$ tends to zero provided 
\begin{equation}
\label{vera}
\vsi\leq\frac{1}{c_3 2^{b+\tau+1}}\ .
\end{equation}
Then the inequalities \eqref{an.5} ensure that the assumptions
\eqref{ite.101} and \eqref{ite.102} of Lemma \ref{iterative} hold. By
taking the limit $l\to\infty$ one gets the result.\qed

\begin{remark}
\label{rem.1.r}
Consider a Hamiltonian of the form 
\begin{equation}
\label{rem.1.p}
H'=H^{(0)}+P'=A^{(0)}+P^{(0)}+P'\ ,
\end{equation}
then the transformation $U\psi'=\psi$ transform it into 
\begin{equation}
\label{rem.2.p}
A^{(\infty)}+U^{-1}P'U\ ,
\end{equation}
and by the estimate \eqref{stitrakam}, one has
\begin{equation}
\label{rem.3.p}
\norma{U^{-1}P'U-P'}_{\vartheta r}\leq 2C_U\vsi r^{b-(n+\tau)}\ ,
\end{equation}
and a similar estimate for the Lipschitz norm.
\end{remark}

\section{KAM with finite smoothness (end of the proof of Theorem
  \ref{m.1})}\label{finite}  

First we define the standard $C^{\ell}$ (H\"older) norms of functions
on $\T^n$ (we use here a definition slightly different from that used
for Whitney smooth functions in order to use tools developed in
\cite{S04}).

Let $0<\mu<1$, and let $F$ be a H\"older function from $\T^n$ to
$B(\cH^{s-\kappa},\cH^s)$ with H\"older exponent $\mu$, then we put
\begin{align}
\label{nmu}
\left|
F\right|_{C^\mu}:=\sup_{0<|\phi-\phi'|<1}\frac{\norma{F(\phi)-F(\phi')}}{|\phi-\phi'|^\mu}+\sup_{\phi\in\T^n}\norma
f\ ,
\\
\label{nmu.1}
\left| F\right|_{C^\ell}:=\sum_{|\alpha|\leq \ell} \left|
\partial^\alpha F\right|_{C^\mu}\ ,\quad \mu:=\ell-\inte \ell \ .
\end{align}

In order to extend a $C^\ell$ function to a complex neighborhood of
$\T^n$ we will use the following polynomials
\begin{equation}
\label{nmu.2}
P_{F,\ell}(\phi,\theta):=\sum_{|\alpha|\leq \ell} \frac{1}{\alpha !}
\partial^\alpha F(\phi) \theta^\alpha\ ,
\end{equation}
and remark that 
\begin{equation}
\label{nmu.3}
\sup_{{|\phi|\in\T^n\atop |\theta|\leq r}}
\norma{P_{F,\ell}(\phi,\theta) }\leq \left| F\right|_{C^\ell} 
\end{equation}
for any $0<r\leq 1$. Then the following smoothing Lemma (from
\cite{S04}) holds.
\begin{lemma}
\label{S3}
(Lemma 3 of \cite{S04}) There is a family of convolution operators
\begin{equation}
\label{S3.1}
{S_r}f(\phi)=\frac{1}{r^n}\int_{\R^n}K\left(\frac{\phi-\phi'}{r}\right)
f(\phi') d^n\phi'\ ,\quad 0<r\leq 1\ ,
\end{equation} 
from $C^0(\R^n)$ into the space of entire analytic functions on $\C^n$
with the following property. For any $\ell\geq 0$ there exists a
constant $c(\ell,n)>0$ such that, for every $\phi\in\C^n$, we have 
$$
|\Im \phi|\leq r\ \Longrightarrow \ \left|\partial^\alpha
S_rf(\phi)-P_{\partial^\alpha f,\inte{\ell}-|\alpha|}(\Re\phi;\ir \Im
\phi)\right| \leq c\left|f\right|_{C^\ell}
r^{\ell-\left|\alpha\right|}\ .
$$ 
Moreover, in $f$ is periodic in $\phi$ then $S_rf$ is periodic in $\Re
\phi$, and $S_rf$ is real valued whenever $f$ is real valued. The
result holds also for functions with values in Banach spaces.
\end{lemma}

A converse of this Lemma is given by

\begin{lemma}
\label{lemma4}
(Lemma 4 of \cite{S04}). Let $\ell\geq 0$ be real, and let $n$ be a positive
integer. Then there exists a constant $c=c(\ell,n)>0$ with the
following property. If $f:\R^n\to\R$ is the limit of a sequence of
 functions $f_\nu(\phi)$ real analytic in the strips $\left|\Im
 \phi\right|\leq r_\nu:=2^{-\nu}r_0$, with $0<r_0\leq  1$ and 
$$
f_0=0\ ,\quad \left|f_{\nu}(\phi)-f_{\nu-1}(\phi)\right|\leq A r_\nu^{\ell}
$$
for $\nu\geq 1$ and $\left|\Im \phi\right|\leq r_\nu$, then $f\in
C^{m}(\R^n)$ for every $m\leq \ell$ which is not integer and moreover
\begin{equation}
\label{smo.4}
|f|_{C^m}\leq\frac{cA}{\mu(1-\mu)}r_0^{\ell-m}\ ,\quad 0<\mu:=m-\inte
m<1 \ .
\end{equation}
\end{lemma}

The proof of the KAM theorem with finite smoothness is based on the
repeated application of Theorem \ref{analyticKAM} to the Hamiltonian
\eqref{H0R}. To describe the procedure we first fix the parameters that
we will use.

Fix a positive $m$ (which will control the smoothness of the
reduction transformation) and some $\ell>m-b$.
Define $\gamma=\ell-m-b$, $r_0:=\epsilon^{1/\ell}$, $r_\nu:=r_0
2^{-\nu}$ and $R^{(\nu)}:=\epsilon S_{r_\nu}R_0$ with $S_r$ the
smoothing operator of Lemma \ref{S3} and $R_0$ the perturbation
defined in \eqref{H0R} (remark that we inserted
$\epsilon$ in the definition of $R^{(\nu)}$). 

The scheme of the iteration is the following one: first we construct
the unitary transformation $U^{(0)}$ transforming $A^{(0)}+ R^{(0)}$
into a time independent diagonal operator $\tilde A^{(1)}$. Then we
use $U^{(0)}$ to transform $ A^{(0)}+ R^{(1)} $. According to Remark
\ref{rem.1.r} it transform such a system into $ \tilde A^{(1)}+ [
  U^{(1)}]^{-1} (R^{(1)}-R^{(0)}) U^{(1)} $, which is a smaller
perturbation of a time independent system. After $\nu-1$ steps we have
thus constructed a unitary transformation $\Phi^{(\nu-1)}$ which
transform $A^{(0)}+ R^{(\nu-1)}$ into $\tilde
A^{(\nu-1)}=\diag(\tilde\lambda_j^{(\nu-1)})$. Use now
$\Phi^{(\nu-1)}$ to transform $A^{(0)}+ R^{(\nu)}$. One gets the
system 
\begin{equation}
\label{pnu.20}
\tilde A^{(\nu-1)}+ [
  \Phi^{(\nu-1)}]^{-1} (R^{(\nu)}-R^{(\nu-1)}) \Phi^{(\nu-1)} \ ,
\end{equation}
to which we apply Theorem \ref{analyticKAM} again. Then one
has to check the assumptions of such a theorem and to add estimates
showing that the procedure converges.

\begin{theorem}
\label{cr}
Consider the quantum Hamiltonian \eqref{H0R}, defined and Lipschitz on
a set $\Omega^{(0)}_{\gamma^{(0)}}$ s.t.
\begin{equation}
\label{an.1.1}
\left|\Omega-\Omega^{(0)}_{\gamma^{(0)}}\right|\leq \Upsilon^{(0)}
\end{equation}
with some positive $\Upsilon^{(0)}\simeq [\gamma^{(0)}]^a$. Let  $\cK_0,$
$\cK_1,$ $\Gamma,$ $\tau,$ be the constants fixed in
Theorem \ref{analyticKAM}. Fix $m$ in such a way that $m+n+\tau$ is
not an integer and let $\ell>m+b$. 
Assume that $R_0\in Lip_1\left(\Omega^{(0)}_\gamma
;C^\ell(\T^n;B(\cH^{s-\kappa},\cH^s))\right)$ (in the following
$C^\ell$ Lipschitz, for short) and let $M$ be a constant such that 
\begin{equation}
\label{sti.1.11}
\norma{R_0}_{C^\ell}\leq M\ ,\quad \normaL{R_0}_{C^\ell}\leq M\ .
\end{equation}
Then there exist positive constants
$\epsilon_*,$ $ C_1$, $C'_{\Omega}$ s.t.,
if $|\epsilon|<\epsilon_*$ and the eigenvalues $\lambda_j^{(0)}$
fulfill \eqref{an.2}-\eqref{an.4}
with constants s.t 
\begin{align}
\label{an.501}
\cK'_{0\lambda}:=\cK^{(0)}_{0\lambda}-2C_0\epsilon^{\frac{m+b}{\ell}}
>\cK_0\ ,\quad \cK'_{1\lambda}:=\cK^{(0)}_{0\lambda}+2C_0
\epsilon^{\frac{b+\tau+1+m}{\ell}} <\cK_1\ , \quad
\gamma':=\gamma^{(0)}-C_{\Gamma} C_1\epsilon^{\frac{b+m}{\ell}d_2}\geq
\Gamma\ ,
\end{align}
then then there exists a measurable set $\widetilde
\Omega^{(\infty)}_{\gamma'}$ 
and a $C^{m+n+\tau}$ Lipschitz map $\Phi^{(\infty)}$ defined on
$\widetilde \Omega^{(\infty)}_{\gamma'}$, with $\Phi^{(\infty)}(\phi,\omega)$
unitary as a map on $L^2$, s.t.  the transformation
$\Phi^{(\infty)}(\omega t,\omega)\psi'=\psi $ transforms the system
\eqref{H0R} into
\begin{equation}
\label{indep.1}
\ir\dot \psi'=\widetilde A^{(\infty)}\psi'\ ,\quad
\widetilde A^{(\infty)}:=\diag(\widetilde \lambda^{(\infty)}_j)\ .
\end{equation}
Furthermore the eigenvalues $\widetilde \lambda^{(\infty)}_j$ fulfill the
estimates \eqref{an.2}-\eqref{an.4} with the new constants defined in
\eqref{an.501} and one has
\begin{align}
\label{an.510}
|\Omega^{(0)}_{\gamma^{(0)}}-\widetilde\Omega^{(\infty)}_{\gamma'}|\leq
C'_{\Omega}\gamma' \ ,
\\
\label{an.511}
\norma{\uno-\Phi^{(\infty)}}_{C^{m'} }\leq 2 C_U \epsilon^{\frac{m
    +n+\tau - m'}{\ell}}\ ,
\end{align}
for any $m'<m+n+\tau$ s.t. $\frac{m+n+\tau-m'}{\ell}$ is not an
integer.  
\end{theorem}

\proof First remark that (by Lemma \ref{S3})
\begin{align}
\label{resti}
\sup_{\phi\in\T^n_{r_\nu}}\norma{R^{(\nu)}(\phi)-\epsilon P_{R_0;\ell}(\Re\phi)\ir
  \Im\phi}\leq c_1\epsilon Mr_\nu^{\ell}\ ,\quad
\sup_{\phi\in\T^n_{r_\nu}}\normaL{R^{(\nu)}(\phi)-P_{R_0,\ell}(\Re\phi)\ir
  \Im\phi}\leq c_1\epsilon Mr_\nu^{\ell}\ .
\end{align}
So that
\begin{equation}
\label{stip0}
\norma{R^{(0)}}_{r_0}\leq c_1\epsilon Mr_0^\ell+\epsilon M\leq 2\epsilon
M= 2M r_0^{\ell}=2M r_0^\gamma r_0^m r_0^b\leq  r_0^m r_0^b\ ,
\end{equation} 
provided $2M r_0^\gamma<1$ (which is a smallness assumption on
$\epsilon$).  For the Lipschitz norm an equal estimate holds:
\begin{equation}
\label{stip0L}
r_0^{\tau+1}\normaL{R^{(0)}}_{r_0}\leq c_1\epsilon Mr_0^\ell+\epsilon
M\leq 2\epsilon M= 2M r_0^{\ell}=2M r_0^\gamma r_0^m r_0^b\leq r_0^m
r_0^b\
\end{equation} 
 (of course one also has a better estimate, but we do not need
it). We also have the following estimates
\begin{align}
\label{pnu}
\norma{R^{(\nu)}-R^{(\nu-1)}}_{r_\nu}\leq
\norma{R^{(\nu)}-P_{R_0;\ell}(\Re\phi)\ir
  \Im\phi}_{r_\nu}+\norma{R^{(\nu-1)}-P_{R_0;\ell}(\Re\phi)\ir
  \Im\phi}_{r_{\nu-1}}
\\
\label{pnu.1}
\leq
c_1\epsilon Mr_\nu^{\ell}+c_1\epsilon Mr_\nu^{\ell}2^{2\ell}=
c_1\epsilon Mr_\nu^{\ell}(1+2^{2\ell})\ .
\end{align}
We now show that for any $\nu$ there exists a set $\widetilde
\Omega^{(\nu)}_{\gamma^{(\nu)}}$ and a Lipschitz analytic transformation $U^{(\nu)}$,
defined on it, unitary in $L^2$ such that, if one defines 
\begin{equation}
\label{pnu.9}
\Phi^{(\nu)}:=U^{(0)}...U^{(\nu)}\ ,
\end{equation}
then it transforms
$A^{(0)}+R^{(\nu)}$ into a time independent system
$A^{(\nu)}=\diag(\tilde\lambda_j^{(\nu)})$ with eigenvalues fulfilling
\eqref{an.2}-\eqref{an.4} with constants  
\begin{align}
\label{pnu.4}
\widetilde{\cK}^{(\nu)}_{0\lambda}=\cK_{0\lambda}^{(0)}-C_0r_0^{b+m}
\sum_{j=0}^{\nu} \frac{1}{2^{(b+m)j}} \ ,\quad \widetilde{\cK}^{(\nu)}_{1\lambda}=\cK_{1\lambda}^{(0)}+C_0r_0^{b+m+\tau+1}
\sum_{j=0}^{\nu} \frac{1}{2^{(b+m+\tau+1)j}}  
\\
\label{pnu.44}
\widetilde\gamma^{(\nu)}=\gamma^{(0)}-C_\Gamma r_0^{(b+m)d_2}\sum_{j=0}^{\nu}
\frac{1}{2^{(b+m)j}}\ ,
\end{align}
and furthermore the following estiamtes hold
\begin{equation}
\label{pnu.2}
\norma{\Phi^{(\nu)}-\uno}_{r_{\nu+1}}\leq
2C_Ur_0^{b_1}\sum_{j=0}^{\nu}\frac{1}{2^{b_1j}}<4C_Ur_0^{b_1}\ ,\quad
r^{\tau+1}\normaL{\Phi^{(\nu)}-\uno}_{r_{\nu+1}}\leq
2C_Ur_0^{b_1}\sum_{j=0}^{\nu}\frac{1}{2^{b_1j}}\ ,
\end{equation}
with $b_1=m+n+\tau$.

Consider the case $\nu=0$. We apply Theorem \ref{analyticKAM} with
$r:=r_0$ and $\vsi:= r_0^m$ to
$A^{(0)}+R^{(0)}$. This is possible since the assumptions on the
eigenvalues are verified by \eqref{an.501}. Then, by
\eqref{stitrakam} and \eqref{an.21}-\eqref{an.51} the equations
\eqref{pnu.2} and \eqref{pnu.4} hold with $\nu=0$. 

Assume now that the result is true for $\nu-1$. Then, as anticipated
above, the transformation $\Phi^{(\nu-1)}$ transform
$A^{(0)}+R^{(\nu)}$ into the system \eqref{pnu.20}, to which we apply
Theorem \ref{analyticKAM}. To this end remark that the assumptions on
the eigenvalues are satisfied by the iterative assumption. We just
have to add an estimate of the new perturbation. In view of the
iterative estimate \eqref{pnu.2} and of \eqref{pnu.1}, one has
$$
\norma{[
  \Phi^{(\nu-1)}]^{-1} (R^{(\nu)}-R^{(\nu-1)})
  \Phi^{(\nu-1)}}_{r_{\nu}}\leq c_1\epsilon
Mr_\nu^\ell(1+2^{2\ell})4\leq r_\nu^mr_\nu^b\ ,
$$ provided $4C_Ur_0^{b_1}<1$ and $c_1\epsilon
Mr_\nu^\gamma(1+2^{2\ell})4<1$, which are smallness assumptions on
$\epsilon$. A similar estimate holds for the Lipschitz norm. 

Applying Theorem \ref{analyticKAM} one gets the transformation
$U^{(\nu)}$ that we need. In
particular the iterative estimate of $\Phi^{(\nu)}$ follows from
\eqref{stitrakam}. 

We have now to show that the sequence of transformations
$\Phi^{(\nu)}$ converges. To this end we apply Lemma \ref{lemma4} to
the sequence $f^{(\nu)}:=\Phi^{(\nu)}-\uno$. Defining $f^{(-1)}:=\uno$
the initial step is fulfilled and one has 
\begin{equation}
\label{pnu.441}
\norma{f^{(\nu)}-f^{(\nu-1)}}_{r_{\nu}}=\norma{\left[U^{(\nu)}-\uno\right]\Phi^{(\nu-1)}}_{r_{\nu}}\leq
2C_Ur_{\nu}^{b_1}\ , 
\end{equation}
which implies the thesis. \qed
 
\noindent{\it End of the proof of Theorem \ref{m.1}} In order to
conclude the proof of Theorem \ref{m.1} one has to show that the
measure of the set of the allowed frequencies becomes full as
$\epsilon\to0$. To this end we remark that the statement implies the
fact that (once all the other parameters are fixed) for any $\Gamma$
there exists $\epsilon_*(\Gamma)>0$ s.t. for smaller $\epsilon$ the
Theorem holds. Denote $d_4:=\frac{\ell}{(b+m)d_2}$. For given
$\Gamma$, define $\gamma^{(0)}:=3\Gamma/2$, then Theorem \ref{cr}
applies provided $\epsilon<\min\left\{C\Gamma^{d_4};\epsilon_*(\Gamma)
\right\}$ with a suitable $C$. Let $\Gamma(\epsilon)$ be the smallest
$\Gamma$ s.t.
$$
2\epsilon=\min\left\{C\Gamma^{d_4};\epsilon_*(\Gamma) \right\}\ ;
$$ we claim that $\Gamma(\epsilon)$ goes to zero. Indeed, assume by
contradiction that this is false, then it means that
$\epsilon_*(\Gamma(\epsilon))>0$ (strictly) for all $\epsilon>0$, but
this contradicts Theorem \ref{cr}. \qed

\section{Proof of Theorem \ref{quasi.l}}\label{p.l}

\noindent {\it Proof of Theorem \ref{quasi.l} in the case $l>1$.} 
Consider the case of $h=h_0+\epsilon W_{2l}$. The lower order
corrections will be added after a first set of transformations.

We start by transforming $h$ using the transformation generated by 
\begin{equation}
\label{p.l.1}
\chi_1:=\frac{b_1(\omega t)x^{l+1}}{l+1}\ ,\quad b_1(\omega t
):=\frac{a_2(\omega t )}{2(1+\epsilon a_1(\omega t))}\ .
\end{equation}
It is easy to see that the flow it generates is
\begin{equation}
\label{p.l.2}
\Phi^\epsilon_{\chi_1}(x,\xi)=(x,\xi-\epsilon b_1(\omega t )x^l)\ ,
\end{equation}
so that, by explicit computation
\begin{equation}
\label{p.l.3}
h^{(1)}:=h\circ \Phi^\epsilon_{\chi_1}= (1+\epsilon a_1(\omega
t))\xi^2+(1+\epsilon c_1(\omega t))x^{2l}\ ,
\end{equation}
with
$$
c_1=a_3+\epsilon b_1^2+\epsilon^2a_1b_1^2-\epsilon a_2b_1\ .
$$
One also has
\begin{equation}
\label{p.l.5}
y^{(1)}_x:=\int_0^\epsilon \dot \chi_1\circ
\Phi_{\chi_1}^{\epsilon-\epsilon_1}d\epsilon_1 +\epsilon\cS^{-(l+1)}= \epsilon
\dot \chi_1+\epsilon\cS^{-(l+1)}\ .
\end{equation}
Remark also that (again by explicit computation)
\begin{equation}
\label{flow.1}
f\circ \Phi^\epsilon_{\chi_1}\in \cS^m\ ,\quad {\rm whenever}\quad
f\in \cS^m\ ,
\end{equation}
so that, by Theorem \ref{7.1}, Eq.\eqref{eg.3} holds. 

In conclusion one has that the transformed Hamiltonian has the form
\begin{equation}
\label{p.l.8}
h^{(1)}-\epsilon\dot \chi_1+\epsilon\cS^{-2}\ .
\end{equation}
We now make a new transformation using 
\begin{equation}
\label{p.l.9}
\chi_2:=b_2(\omega t)x\xi\ ,\quad b_2(\omega t):=\frac{1}{4(l+1)\epsilon}\ln\left(\frac{1+\epsilon
  a_1}{1+\epsilon c_1}\right)\ ,
\end{equation}
(the $\epsilon$ in $b_2$ only plays the role of a parameter)
whose flow is given by
$$
\phi^\epsilon_{\chi_2}(x,\xi)=(e^{b_2\epsilon}x,e^{-b_2\epsilon}\xi)\ .
$$
thus Eq. \eqref{flow.1} holds. One has
\begin{equation}
\label{p.l.10}
h^{(2)}:=h^{(1)}\circ\phi^\epsilon_{\chi_2}=c_2(\omega
t)(\xi^2+x^{2l})\ ,\quad c_2=(1+\epsilon
c_1)^{\frac{1}{2(l+1)}}(1+\epsilon a_1)^{\frac{2l+1}{2l+2}} \ ,
\end{equation}
and 
$$
y^{(2)}_x=-\epsilon\dot \chi_2+\epsilon \cS^{-(l+1)}\ .
$$
Thus, after this couple of transformations $h_0+\epsilon W_{2l}$ is
transformed to
$$
c_2(\omega t)h_0-\epsilon p_{l+1}+\epsilon\cS^{-2}\ ,
$$
with
$$
p_{l+1}:=\dot \chi_1\circ \phi^\epsilon_{\chi_2}+\dot \chi_2\ ,
$$ which is quasihomogeneous of degree $l+1$. The idea (following
\cite{BBM14}) is now to get rid of the time dependence of the main term
by reparametrizing time, i.e. to pass to a new time $\tau$ such that 
\begin{equation}
\label{p.l.11}
\frac{ d \tau}{dt}=c_2(\omega t)\ .
\end{equation}
First we show that \eqref{p.l.11} defines a good reparametrization of
time. Indeed, by making a Fourier expansion of $c_2$:
$$
c_2(\phi)=(1+\epsilon c_{20})+\epsilon \sum_{k\not =0} c_{2k} e^{\ir
  k\cdot \phi}\ ,
$$
one has
\begin{equation}
\label{p.l.13}
\tau(t)=(1+\epsilon c_{20})t+\epsilon  \sum_{k\not =0}
\frac{c_{2k}}{\ir \omega\cdot k} e^{\ir
  k\cdot\omega t} =:(1+\epsilon c_{20})t+\epsilon \tau_1(\omega t)\ ,
\end{equation}
which is well defined and $C^\infty$ on $\Omega_{0\gamma}$. Then one
can use the implicit function theorem in order to show that the
inverse of the transformation \eqref{p.l.13} has the form 
\begin{equation}
\label{p.l.14}
t(\tau)=a\tau - \epsilon t_1(\omega a \tau)\ , \quad a:= (1+\epsilon
c_{20})^{-1}\ ,
\end{equation}
and $t_1$ defined and smooth on $\T^n$. Precisely this is obtained by
applying the implicit function to the equation (that defines $t_1$)
$$
G(\epsilon,t_1):=\tau_1(\phi -\epsilon \omega t_1(\phi)
)-t_1(\phi)=0\ ,
$$
where $G:\R\times C^K_0(\T^n)\to C^K_0(\T^n)$, with an arbitrary
$K$ and the index $0$ means ``with zero average''. 

After the introduction of the new time the system is reduced to the
quantization of 
$$
h_0+\epsilon p'_{l+1}+\cS^{-2}\ ,
$$ 
with $p'_{l+1}:=p_{l+1}/c_2$ (and the frequencies are now
$\omega'=a\omega$). 

We have now to eliminate $p'_{l+1}$. To this end we proceed as
explained in Sect. \ref{smoa}, i.e. we solve eq. \eqref{12} with
$p=p'_{l+1}$, thus getting a $\chi_3\in S^2$ which conjugates the
Hamiltonian to
$$
h_0+\epsilon \langle p'_{l+1}\rangle (h_0,\omega' \tau)\eta(h_0)
+S^2\ .
$$
The last step is achieved by removing the time dependence from $
\langle p'_{l+1}\rangle\eta$. To this end we look for a $\chi_4$
solving \eqref{12.p} with $p=\langle p'_{l+1}\rangle\eta$. The main
remark is that the function $\chi_4\in S^{l+1}$ turns out to be
quasihomogeneous (in the region $E>2$), thus it is easy to see that it
has the property that $f\circ\Phi^{\epsilon}_{\chi_4}\in S^m$ whenever
$f\in S^m$, and therefore eq. \eqref{eg.3} holds. Using such a
$\chi_4$ one conjugates the Hamiltonian to
$$
h_0+\epsilon \overline{\langle p'_{l+1}\rangle}+\epsilon S^2\ .
$$
At this point we can add the lower order corrections $W$ and apply
Theorem \ref{m.1} getting the result. \qed
\vskip 10 pt

\noindent {\it Proof of Theorem \ref{quasi.l} in the case $l=1$.}  The
proof is a simple KAM type theorem in which, working at the level of
symbols, one eliminates iteratively the time dependence from the
Hamiltonian. The key remark is that, if $\chi$ is quadratic, then
given a symbol $h_0+\epsilon p$, one has that the
symbol of $T_{\epsilon X}(H_0+\epsilon P)$ is {\it exactly} 
\begin{align}
\label{quasi.1}
(h_0+\epsilon p)\circ \Phi^\epsilon_{\chi}-y_x\simeq
h_0+\epsilon\left\{ h_0;\chi\right\}-\epsilon\dot\chi\ ,
\\
\label{quasi.2}
y_x=\int_0^\epsilon\dot\chi\circ\Phi^{\epsilon-\epsilon_1}_\chi d\epsilon_1\ .
\end{align}
So, in order to establish the recursion one determines $\chi$ by
solving the homological equation \eqref{12.1} (with $p$ in place of
$W$) and uses it in order square the order of the time dependent part
of the symbol. 

Then one has to add estimates and to prove an iterative Lemma which
allows to establish the convergence of the procedure. We remark that
such an iterative Lemma is actually a simple 2-dimensional version of
Lemma \ref{iterative}. For this reason we omit the details of the
proof. \qed

\appendix

\section{Some technical lemmas}\label{lemms}

We start with a couple of results which apply to operators depending
in a $C^\infty$ way on the angles. They are used in sect. \ref{prepa}.

\begin{remark}
\label{A.commu}
For a  rough estimate of the commutator of two operators remark that,
having fixed a set $\tilde \Omega$ and indexes $\rho,s,\delta$, then
there exists a constant $C$ which depends on all these indexes s.t.
\begin{equation}
\label{A.commu.1}
\norma{\left[X,F\right]}_{Lip_\rho\left(\tilde
  \Omega;B(\cH^{s-\delta};\cH^s)\right)}\leq C \norma X_{Lip_\rho\left(\tilde
  \Omega;B(\cH^{s-\delta};\cH^s)\right)} \norma F_{Lip_\rho\left(\tilde
  \Omega;B(\cH^{s-\delta};\cH^s)\right)} \ .
\end{equation}
\end{remark}

Exploiting such a Remark it is immediate to get the following result
whose proof is obtained just by estimating each term of the series
defining the quantum Lie transform and summing up the series.

\begin{lemma}
\label{qliest}
Let $F$ and $X$ be two operators belonging to $ Lip_\rho\left(\tilde
  \Omega;B(\cH^{s-\delta};\cH^s)\right)$, then $Lie_{\epsilon X}F$ also belongs
  to such a space and there exists a constant
  $C$ which depends only on the indexes of the norm, such that
\begin{equation}
\label{qliest.1}
\norma{Lie_{\epsilon X}F-F}\leq C \norma X\norma F\ .
\end{equation}
The norm is the norm in the above space.
\end{lemma}

We prove now some general properties of sequences $\lambda_j$, having a
behaviors of that of the eigenvalues of the operators that we meet in
the main part of the text.

\begin{lemma}
\label{meas.1}
Assume that 
\begin{equation}
\label{meas.2}
\left|\lambda_i(\omega)-\lambda_j(\omega)\right|\geq\cK_{0\lambda}
\left|i^d-j^d 
\right|\ ,\quad i\not=j\ ,
\end{equation}
and $\alpha\leq \cK_{0\lambda}/2$, then
\begin{equation}
\label{notempty}
\cR_{ijk}(\lambda,\alpha)\not=\emptyset\ \Longrightarrow
\ \left|k\right| \geq \frac{\cK_{0\lambda}}{4} \left|i^d-j^d 
\right|\ .
\end{equation}
\end{lemma}
\proof Since $\cR_{ijk}(\lambda,\alpha)\not=\emptyset$ one has
$$
2\left| k\right|\geq \left|\omega\right|_{\ell^\infty} |k|\geq
 \left|\lambda_i-\lambda_j\right|-\alpha \left|i^d-j^d\right|\geq
 (\cK_{0\lambda}-\alpha) \left|i^d-j^d\right|\ .
$$
\qed

\begin{lemma}
\label{meas.4}
Assume \eqref{meas.2} and
\begin{equation}
\label{meas.22}
\left|\frac{\Delta(\lambda_i-\lambda_j)}{\Delta\omega}\right|\leq
\cK_{1\lambda}\idjd \ ,\quad \forall \omega\in \tilde\Omega\ ,
\end{equation}
with $\cK_{1\lambda}\leq \cK_{0\lambda}/8$, then one has
\begin{equation}
\label{meas.23}
\left|\cR_{ijk}(\lambda,\alpha)\right|\leq
\frac{4\alpha}{\cK_{0\lambda}}n^{(n-1)/2}\ .
\end{equation}
\end{lemma}
\proof Assume that $\cR_{ijk}$ is not empty, so that 
\eqref{notempty} holds.  Let $\omega\in\cR_{ijk}$; choose a vector
$v\in\left\{-1,1\right\}^n$ such that $k\cdot v=|k|$ and write
$\omega=rv+w$ with $w\in v^{\perp}$. We estimate the size by which one
has to move $r$ in order to go outside $\cR$. Let $\omega':=r'v+w$
and compute
\begin{align*}
\left|\Delta s'_{ijk}\right|\geq
|(r-r')||k|-|\Delta(\lambda_i-\lambda_j)|\geq |\Delta\omega|\left( |k|-
\cK_{1\lambda}\idjd\right)\ .
\end{align*} 
So, if such a quantity is larger then $\alpha\idjd$, then $\omega'$ is
outside $\cR_{ijk}$. It follows that  
\begin{align*}
\left|\cR_{ijk}\right|\leq \frac{2\alpha\idjd}{|k|-\cK_{1\lambda}
  \idjd}(diam(\Omega))^{n-1} \leq
\frac{2\alpha\idjd\cK_{0\lambda}}{(\cK_{0\lambda}-4\cK_{1\lambda})|k|
}n^{(n-1)/2} \leq
\frac{4\alpha}{\cK_{0\lambda}}n^{(n-1)/2}\ .
\end{align*}
\qed

\section{Estimates of analytic quantum Lie transform}\label{analiticq}

\begin{remark}
\label{commu}
One has 
\begin{align}
\label{commu.1}
\norma{\ir [X;F]}_{r-\sigma}&\leq 2\norma{X}_{r-\sigma}\norma{F}_{r-\sigma} 
\\
\label{commu.2}
\normaL{\ir [X;F]}_{r-\sigma}&\leq
2\left(\normaL{X}_{r-\delta}\norma{F}_{r-\sigma}+\norma{X}_{r-\sigma}
\normaL{F}_{r-\sigma}\right)\ .
\end{align}
\end{remark}

\begin{lemma}
\label{reLie}
Provided 
\begin{equation}
\label{smaReLie}
\norma X_{r_\sigma}<\frac{\ln 2}{2}\ ,
\end{equation}
one has
\begin{align}
\label{reLi1.1}
\norma{e^{\ir X}Fe^{-\ir X}-F}_{r-\sigma}&\leq
4\norma{X}_{r-\sigma}\norma F_{r}\ ,
\\
\label{reLie.2}
\normaL{e^{\ir X}Fe^{-\ir X}-F}_{r-\sigma}&\leq
4\norma{X}_{r-\sigma}\normaL F_{r}+2\normaL{X}_{r-\sigma}\norma F_{r}\ ,
\end{align}
\end{lemma}
\proof From the recursive formula \eqref{serqlie} and remark
\ref{commu}, one immediately
gets 
\begin{equation}
\label{fkest}
\norma{F_k}_{r-\sigma}\leq (2\norma{X}_{r-\sigma})^k\norma{F}_{r-\sigma}\ ,
\end{equation}
from which 
$$
\norma{e^{\ir X}Fe^{-\ir X}-F}_{r-\sigma}\leq\sum_{k\geq1}{\frac{(
    2\norma{X}_{r-\sigma})^k }{k!}}\norma{F}_r\leq (e^{
    2\norma{X}_{r-\sigma}} -1)\norma{F}_r\ ,
$$
which, under the assumption \eqref{smaReLie}, is smaller then the
r.h.s. of \eqref{reLi1.1}.

We come to \eqref{reLie.2}. From \eqref{commu.2} one gets
\begin{align}
\label{Fklie}
\normaL{F_k}_{r-\sigma}&\leq2\left(\norma{X}_{r-\sigma}\normaL{F_{k-1}}_{r-\sigma}+
\normaL{X}_{r-\sigma}\norma{F_{k-1}}_{r-\sigma}\right)
\\
& \leq
2\norma{X_{r-\sigma}}\normaL{F_{k-1}}_{r-\sigma}+\normaL{X}_{r-\sigma}(2\norma{X}_{r-\sigma})^{k-1}\norma{F}_r
\ .
\end{align}
To write the formulae we need in a simpler way denote
\begin{equation}
\label{notazioni}
\lambda:=2\norma{X}_{r-\sigma}\ ,\quad
\mu:=\normaL{X}_{r-\sigma}\ , b:=\norma{F}_r\ ,
\end{equation}
and look for a sequence $a_k$ such that $a_k\geq
\normaL{F_k}_{r-\sigma}$, then such a sequence can be defined by
\begin{equation}
\label{ak}
a_k=\lambda a_{k-1}+\mu\lambda^{k-1}b\ ,
\end{equation}
which is easily solved by the discrete equivalent of Duhamel formula,
which is actually obtained by making the substitution $a_k=\lambda^k
c_k$, so that $c_k$ satisfies
$$
c_k=c_{k-1}+\frac{\mu}{\lambda}b\ ,
$$
which gives
$$
c_k=c_0+k\frac{\mu}{\lambda}b\ ,\quad
a_k=\lambda^ka_0+\lambda^{k-1}k\mu b\ ,
$$
Now, the l.h.s. of \eqref{reLie.2} is estimated by
$$
\sum_{k\geq1}\frac{a_k}{k!}=(e^{\lambda}-1)a_0+\mu be^{\lambda}\ ,
$$
which, again under \eqref{smaReLie} gives the result.\qed

Let $X$ be the solution of Eq. \eqref{qhom}
then, if $A$ is not bounded its Lie transform with $X$ has good
properties. Indeed the following Lemma holds:
\begin{lemma}
\label{estia}
One has
\begin{align}
\label{estia.1}
&\norma{e^{\ir X}Ae^{-\ir X}-A-\ir\left[ X;A\right]}_{r-2\sigma}\leq
4\norma{X}_{r-\sigma}\left[\frac{1}{\sigma}\norma{X}_{r-\sigma}+2\norma{P}_{r-2\sigma}\right]\ ,
\\
\nonumber
&\normaL{e^{\ir X}Ae^{-\ir X}-A-\ir\left[ X;A\right]}_{r-2\sigma}
\\
\label{estia.2}
&\null\qquad\leq
\frac{8}{\sigma}\norma{X}_{r-\sigma}\normaL{X}_{r-\sigma}
+\frac{4}{\sigma}\norma{X}_{r-\sigma}^2
+8\norma{X}_{r-\sigma}\normaL{P}_{r-2\sigma}+8 
\normaL{X}_{r-\sigma}\norma{P}_{r-2\sigma}\ .
\end{align}
\end{lemma}
\proof Just remark that the recursion defining $A_k$ can be generated
starting from 
$$
A_1=-\ir\left[A;X\right]=\dot X-P+[P]\ ,
$$
which allows to start with the estimate
\begin{equation}
\label{estia.3}
\norma{A_1}_{r-2\sigma}\leq
\frac{\left|\omega|\right|}{\sigma}\norma{X}_{r-\sigma}+2
\norma{P}_{r-2\sigma}=:b   \ ;
\end{equation}
from this one gets $\norma{A_k}_{r-2\sigma}\leq \lambda^{k-1}b$ with
$\lambda$ defined by \eqref{notazioni}, which gives that the l.h.s. of
\eqref{estia.1} is estimated by
$$
\sum_{k\geq2}\frac{\lambda^{k-1}b}{k!}=\frac{b}{\lambda}(e^\lambda-1-\lambda)\ ,
$$
which gives the wanted estimate.

To estimate the Lipschitz norm we proceed as above: let $a_k$ be a
sequence estimating the Lipschitz norm of $A_k$, then we have
$$
a_k=2b\mu\lambda^{k-2}+\lambda
a_{k-1}=\frac{2b\mu}{\lambda}\lambda^{k-1}+\lambda b_{k-1} \ .
$$
Proceeding again by discrete Duhamel formula ($a_k=\lambda^{k-1}c_k$),
one gets 
$$
c_k=\frac{2b\mu}{\lambda}+c_{k-1}\ ,
$$
which gives
$$
c_k=\frac{2b\mu}{\lambda}(k-1)+a_1\ \,\quad
a_k=\lambda^{k-1}\left(\frac{2b\mu}{\lambda}(k-1)+a_1\right) 
$$
It follows that
\begin{align*}
\sum_{k\geq2}\frac{a_k}{k!}=\sum_{k\geq2}\lambda^{k-1}\frac{2b\mu}{\lambda
k!}-\sum_{k\geq2}\left(\frac{2b\mu}{\lambda
k!}\frac{\lambda^{k-1}}{\lambda}-\frac{a_1\lambda^{k-1}}{k!}
\right)
\\
=\frac{2b\mu}{\lambda}(e^\lambda-1)+\left(a_1-\frac{2b\mu}{\lambda}\right)\frac{e^\lambda-1-\lambda}{\lambda}
\\
=\frac{2b\mu}{\lambda}\left(\frac{\lambda
  e^\lambda-e^\lambda+1}{\lambda}\right)
+a_1\frac{e^\lambda-1-\lambda}{\lambda} \leq \frac{2b\mu}{\lambda}
\lambda^2e^\lambda +2a_1\lambda\leq 4b\mu+2a_1\lambda
\end{align*}
from which, taking
$$
a_1:=\frac{1}{\sigma}\normaL{X}_{r-\sigma}+ 
\frac{1}{\sigma}\norma{X}_{r-\sigma}+2\normaL{P}_{r-2\sigma}
$$ the thesis follows\qed

Concerning $Y_X$ we have the following Lemma:

\begin{lemma}
\label{lemmay}
Let $Y_X$ be defined by \eqref{yx} with $\epsilon=1$, then we have
the following estimates
\begin{align}
\label{estiy}
\norma{Y_X+\dot X}_{r-2\sigma}&\leq\frac{4}{\sigma}\norma{X}_{r-\sigma}^2
\\
\label{estiy1}
\normaL{Y_X+\dot
  X}_{r-2\sigma}&\leq\frac{6}{\sigma}\norma{X}_{r-\sigma} \normaL{X}_{r-\sigma}
\end{align}
\end{lemma}\proof Define $Y_k$ by 
$$
Y_0=\dot X\ ,\quad Y_k:=\ir[X;Y_{k-1} ]\ ,
$$ then
\begin{align*}
Y_X=-\int_0^1 d\epsilon_1e^{\ir(1-\epsilon_1)X}\dot X
e^{-\ir(1-\epsilon_1)X}=-\int_0^1 d\epsilon_1
\sum_{k\geq0} Y_k\frac{(1-\epsilon_1)^k}{k!} 
\\
=-\sum_{k\geq0}Y_k\frac{1}{k!}(-)\frac{1}{k+1}(1-\epsilon_1)
\big|_{0}^1  =-\sum_{k\geq 0}Y_k\frac{1}{(k+1)!}
\ .
\end{align*}
Thus, following the proof of Lemma \ref{reLie}, one easily gets the
thesis. \qed

\section{A few more lemmas}\label{lemkam}

Finally we add two lemmas on the solution of the quantum homological equation
\eqref{cohomoq} and a lemma allowing to solve superexponential recursions.

\begin{lemma}
\label{pos96} (Lemma A.1 of \cite{pos96})
If $P=(P_{ij})$ is (the matrix of) a bounded linear operator on
$\ell^2$, then also $X=(X_{ij})$ with
\begin{equation}
\label{pos96.1}
X_{ij}=\frac{|P_{ij}|}{|i-j|}\ ,\quad i\not=j
\end{equation}
and $B_{ii}=0$ is a bounded linear operator on $\ell^2$, and
$\norma{X}\leq \pi\norma{P}/\sqrt3$, where the norm is the norm in
$B(\ell^2;\ell^2)$.
\end{lemma}

\begin{corollary}
\label{pos92.c}
Let $P\in B(\cH^{s_1};\cH^{s_2})$ for some $s_1,s_2$ be a bounded
operator with matrix $P_{ij}$. Define $X$ by \eqref{pos96.1}, then
also $X\in B(\cH^{s_1};\cH^{s_2})$ and one has
\begin{equation}
\label{pos96.4}
\norma{X}_{B(\cH^{s_1};\cH^{s_2})}\leq
\frac{\pi}{\sqrt3}\norma{P}_{B(\cH^{s_1};\cH^{s_2})}\ .
\end{equation}
\end{corollary}
\proof  Remark that an operator $P$ belongs to
$B(\cH^{s_1};\cH^{s_2})$ if and only if the operator with matrix
$i^{s_2}P_{ij}j^{-s_1}$ is bounded on $\ell^2$, and apply Lemma
\ref{pos96}. \qed

\begin{lemma}
\label{rico}
For $\nu\geq0$, define
\begin{equation}
\label{rico.1}
\vsi_{\nu+1}:=c_12^{a\nu}\vsi_\nu^2\ ,
\end{equation} 
then one has
\begin{equation}
\label{rico.2}
\vsi_\nu=\frac{1}{c_12^{a\nu}}\left(2^{2a}c_1\vsi_0\right)^{2^{\nu}}\ .
\end{equation}
Assume also $2^{2a}c_1\vsi_0<1$, then one has
\begin{equation}
\label{rico.3}
\sum_{\nu\geq
  k}\vsi_\nu\leq \frac{(2^{2a} c_1\vsi_0)^{2^{k}}2}{c_1
  2^{ak}} \ ,
\end{equation}
for any $b>0$ there exists $C_b$ independent of $c_1$  s.t.
\begin{equation}
\label{rico.4}
\sum_{\nu\geq0}\nu^b\vsi_\nu\leq C_b\vsi_0\ .
\end{equation}
\end{lemma}
\proof Make the substitution $\vsi_\nu=k_1k_2^\nu\delta_\nu$, and
rewrite formula \eqref{rico.1} for the sequence $\delta_\nu$. One gets
$$
k_1k_2^{\nu+1}\delta_{\nu+1}=2^{a\nu}c_1k_1^2k_2^{2\nu}\delta_\nu^2\ ,
$$
which becomes particularly simple taking 
\begin{align*}
k_2^{\nu}=2^{a\nu}k_2^{2\nu}\ \iff\ k_2=2^{-a}\ ,
\\
k_2k_1=c_1k_1^2\ \iff\ k_1=(2^ac_1)^{-1}\ . 
\end{align*}
so that we get 
$$
\delta_{\nu+1}=\delta_\nu^2\ \iff\ \delta_\nu=\delta_1^{2^{\nu-1}}\ .
$$
Substituting back in $\vsi_\nu$ one gets \eqref{rico.2}. To get
\eqref{rico.3} remark that
$$ \sum_{\nu\geq k}\vsi_\nu= \frac{(2^{2a}
  c_1\vsi_0)^{2^{k}}}{c_1 2^{ak}}\sum_{\nu\geq
  k}\frac{1}{2^{a(\nu-k)}}(2^{2a}c_1\vsi_0)^{2^{\nu}-2^{k}}\ .
$$
remark that 
$$
2^{\nu}-2^{k}=2^{k}(2^{\nu-k}-1)\geq 2^{k}(\nu-k)\ ,
$$
so that the above sum is smaller than
$$
\frac{(2^{2a}
  c_1\vsi_0)^{2^{k}}}{c_1 2^{ak}}\sum_{\nu\geq
  k}\frac{1}{2^{a(\nu-k)}}(2^{2a}c_1\vsi_0)^{2^{k}(\nu-k)}\leq \frac{(2^{2a}
  c_1\vsi_0)^{2^{k}}}{c_1 2^{ak}}2\ .
$$
The cases with $b>0$ are estimated in the same way.
\qed

\addcontentsline{toc}{chapter}{Bibliography}


\begin{thebibliography}{DL{\v{S}}V02}

\bibitem[Bam97]{Bam97}
D.~Bambusi.
\newblock Long time stability of some small amplitude solutions in nonlinear
  {S}chr\"odinger equations.
\newblock {\em Comm. Math. Phys.}, 189(1):205--226, 1997.

\bibitem[BBM14]{BBM14}
P.~Baldi, M.~Berti, and R.~Montalto.
\newblock K{AM} for quasi-linear and fully nonlinear forced perturbations of
  {A}iry equation.
\newblock {\em Math. Ann.}, 359(1-2):471--536, 2014.

\bibitem[BG93]{BG93}
D.~Bambusi and A.~Giorgilli.
\newblock Exponential stability of states close to resonance in
  infinite-dimensional {H}amiltonian systems.
\newblock {\em J. Statist. Phys.}, 71(3-4):569--606, 1993.

\bibitem[BG01]{BG01}
D.~Bambusi and S.~Graffi.
\newblock Time quasi-periodic unbounded perturbations of {S}chr\"odinger
  operators and {KAM} methods.
\newblock {\em Comm. Math. Phys.}, 219(2):465--480, 2001.

\bibitem[BGP99]{BGP99}
D.~Bambusi, S.~Graffi, and T.~Paul.
\newblock Normal forms and quantization formulae.
\newblock {\em Comm. Math. Phys.}, 207(1):173--195, 1999.

\bibitem[BM16a]{BM16}
D.~Bambusi and A.~Maspero.
\newblock Freezing of {E}nergy of a {S}oliton in an {E}xternal {P}otential.
\newblock {\em Comm. Math. Phys.}, 344(1):155--191, 2016.

\bibitem[BM16b]{BM16a}
M.~Berti and R.~Montalto.
\newblock Quasi-periodic standing wave solutions of gravity-capillary water
  waves.
\newblock {\em arXiv:1602.02411 [math.AP]}, 2016.

\bibitem[Com87]{C87}
M.~Combescure.
\newblock The quantum stability problem for time-periodic perturbations of the
  harmonic oscillator.
\newblock {\em Ann. Inst. H. Poincar\'e Phys. Th\'eor.}, 47(1):63--83, 1987.

\bibitem[Del14]{D14}
J.-M. Delort.
\newblock Growth of {S}obolev norms for solutions of time dependent
  {S}chr\"odinger operators with harmonic oscillator potential.
\newblock {\em Comm. Partial Differential Equations}, 39(1):1--33, 2014.

\bibitem[DL{\v{S}}V02]{DSV02}
P.~Duclos, O.~Lev, P.~{\v{S}}{\v{t}}ov{\'{\i}}{\v{c}}ek, and M.~Vittot.
\newblock Weakly regular {F}loquet {H}amiltonians with pure point spectrum.
\newblock {\em Rev. Math. Phys.}, 14(6):531--568, 2002.

\bibitem[D{\v{S}}96]{DS96}
P.~Duclos and P.~{\v{S}}{\v{t}}ov{\'{\i}}{\v{c}}ek.
\newblock Floquet {H}amiltonians with pure point spectrum.
\newblock {\em Comm. Math. Phys.}, 177(2):327--347, 1996.

\bibitem[EK09]{EK09}
H.~L. Eliasson and S.~B. Kuksin.
\newblock On reducibility of {S}chr\"odinger equations with quasiperiodic in
  time potentials.
\newblock {\em Comm. Math. Phys.}, 286(1):125--135, 2009.

\bibitem[FGJS04]{FGJS04}
J.~Fr{\"o}hlich, S.~Gustafson, B.~L.~G. Jonsson, and I.~M. Sigal.
\newblock Solitary wave dynamics in an external potential.
\newblock {\em Comm. Math. Phys.}, 250(3):613--642, 2004.

\bibitem[FP15]{FP15}
R.~Feola and M.~Procesi.
\newblock Quasi-periodic solutions for fully nonlinear forced reversible
  {S}chr\"odinger equations.
\newblock {\em J. Differential Equations}, 259(7):3389--3447, 2015.

\bibitem[GP87]{GP87}
S.~Graffi and T.~Paul.
\newblock The {S}chr\"odinger equation and canonical perturbation theory.
\newblock {\em Comm. Math. Phys.}, 108(1):25--40, 1987.

\bibitem[GP16]{GP16}
B.~Gr\'ebert and E.~Paturel.
\newblock On reducibility of quantum harmonic oscillator on ${\R}^d$ with
  quasiperiodic in time potential.
\newblock {\em arXiv:1603.07455 [math.AP]}, 2016.

\bibitem[GT11]{GT11}
B.~Gr{\'e}bert and L.~Thomann.
\newblock K{AM} for the quantum harmonic oscillator.
\newblock {\em Comm. Math. Phys.}, 307(2):383--427, 2011.

\bibitem[GY00]{GY00}
S.~Graffi and K.~Yajima.
\newblock Absolute continuity of the {F}loquet spectrum for a nonlinearly
  forced harmonic oscillator.
\newblock {\em Comm. Math. Phys.}, 215(2):245--250, 2000.

\bibitem[HR82a]{HR82D}
B.~Helffer and D.~Robert.
\newblock Asymptotique des niveaux d'\'energie pour des hamiltoniens \`a un
  degr\'e de libert\'e.
\newblock {\em Duke Math. J.}, 49(4):853--868, 1982.

\bibitem[HR82b]{HR82}
B.~Helffer and D.~Robert.
\newblock Propri\'et\'es asymptotiques du spectre d'op\'erateurs
  pseudodiff\'erentiels sur {${\bf R}^{n}$}.
\newblock {\em Comm. Partial Differential Equations}, 7(7):795--882, 1982.

\bibitem[IPT05]{IPT05}
G.~Iooss, P.~I. Plotnikov, and J.~F. Toland.
\newblock Standing waves on an infinitely deep perfect fluid under gravity.
\newblock {\em Arch. Ration. Mech. Anal.}, 177(3):367--478, 2005.

\bibitem[LY10]{LY10}
J.~Liu and X.~Yuan.
\newblock Spectrum for quantum {D}uffing oscillator and small-divisor equation
  with large-variable coefficient.
\newblock {\em Comm. Pure Appl. Math.}, 63(9):1145--1172, 2010.

\bibitem[Mon14]{Mtesi}
R.~Montalto.
\newblock {KAM} for quasi-linear and fully nonlinear perturbations of {A}iry
  and {K}d{V} equations.
\newblock {\em {P}hd {Thesis}, SISSA - ISAS}, 2014.

\bibitem[MR16]{MR16}
A.~Maspero and D.~Robert.
\newblock On time dependent {S}chr\"odinger equations: global well-posedness
  and growth of {S}obolev norms.
\newblock {\em Preprint}, 2016.

\bibitem[P{\"o}s96]{pos96}
J. P{\"o}schel.
\newblock A {KAM}-theorem for some nonlinear partial differential equations.
\newblock {\em Ann. Scuola Norm. Sup. Pisa Cl. Sci. (4)}, 23(1):119--148, 1996.

\bibitem[PT01]{PT01}
P.~I. Plotnikov and J.~F. Toland.
\newblock Nash-{M}oser theory for standing water waves.
\newblock {\em Arch. Ration. Mech. Anal.}, 159(1):1--83, 2001.

\bibitem[Sal04]{S04}
D.~A. Salamon.
\newblock The {K}olmogorov-{A}rnold-{M}oser theorem.
\newblock {\em Math. Phys. Electron. J.}, 10:Paper 3, 37 pp. (electronic),
  2004.

\bibitem[Ste70]{Stein}
E.~M. Stein.
\newblock {\em Singular integrals and differentiability properties of
  functions}.
\newblock Princeton Mathematical Series, No. 30. Princeton University Press,
  Princeton, N.J., 1970.

\bibitem[Wan08]{W08}
W.-M. Wang.
\newblock Pure point spectrum of the {F}loquet {H}amiltonian for the quantum
  harmonic oscillator under time quasi-periodic perturbations.
\newblock {\em Comm. Math. Phys.}, 277(2):459--496, 2008.

\bibitem[YZ13]{YZ13}
X.~Yuan and K.~Zhang.
\newblock A reduction theorem for time dependent {S}chr\"odinger operator with
  finite differentiable unbounded perturbation.
\newblock {\em J. Math. Phys.}, 54(5):052701, 23, 2013.

\end{thebibliography}

\end{document}